\documentclass{article}
\usepackage{graphicx} 
\usepackage[hidelinks]{hyperref}
\usepackage{authblk}
\usepackage{color}
\usepackage{graphicx,subcaption} 
\usepackage{amsmath}
\usepackage{multirow}
\usepackage{float}
\usepackage{hhline}
\usepackage{amsthm} 
\usepackage{amsfonts}
\usepackage{mathtools}
\numberwithin{equation}{section}

\usepackage[table,xcdraw]{xcolor}
\title{Boundary-Layer-Induced Failure of Standard Physics-Informed Neural Networks: A Legendre Wavelet Collocation Benchmark for Singularly Perturbed Transport Problems}
\author[1]{Suvendu Nayak\thanks{Corresponding author: suvendu.nayak1601@gmail.com}}
\author[2]{Arun Kumar Gupta\thanks{arun27.nit@gmail.com}}
\affil[1,2]{Department of Mathematics and Applied Statistics, School of Applied Sciences, KIIT Deemed to be University, Bhubaneswar, Odisha, India}
\usepackage[paper=a4paper,total={7in,10in}]{geometry}

\begin{document}
		\maketitle
\begin{abstract}
	Boundary layers provide a demanding test for numerical solvers because the solution may remain almost constant over most of the domain while changing rapidly in a narrow region near the boundary. This paper studies a singularly perturbed one-dimensional transport boundary-value problem with increasing Peclet number $(\mathrm{Pe})$. A local Legendre wavelet collocation method (LWM) is compared with a standard soft-boundary physics-informed neural network (PINN) for this benchmark. The wavelet approximation uses locally supported Legendre polynomial basis functions and converts the problem into a square algebraic collocation system with residual, boundary, and interface-continuity equations. Numerical experiments are performed for $\mathrm{Pe}=1,10,100,$ and $1000$. The LWM captures all four cases, with the largest error remaining below $5\times 10^{-3}$. The standard soft-boundary PINN performs well for the mild cases, but fails to resolve the sharp boundary layer for the larger Peclet numbers. The results show that local wavelet collocation is more reliable than the standard soft-boundary PINN for this benchmark, while dense near-boundary evaluation helps reveal errors that may be missed on coarse grids.
\end{abstract}


\section{Introduction}

Singularly perturbed transport problems arise when diffusion is small compared with convective or drift effects. In such problems, the solution may remain nearly flat over most of the domain and then change rapidly inside a narrow region close to a boundary. This narrow transition region is usually called a boundary layer. The numerical difficulty is not only caused by the differential equation itself, but also by the fact that the important variation of the solution is localized in a very small part of the interval. For this reason, singularly perturbed boundary-value problems have long been used as demanding tests for numerical methods \cite{Roos2008,Miller1996,Farrell2000,Kopteva2010}.

A second difficulty is the reliable measurement of error. For boundary-layer problems, a coarse uniform set of evaluation points can give a misleading impression of accuracy because most of the interval may contain little variation. The largest error may occur only near the layer, where the solution changes rapidly. Therefore, robust assessment requires either layer-resolving grids or additional evaluation points near the boundary layer. This issue is well known in the numerical analysis of singularly perturbed convection--diffusion problems and fitted-mesh methods \cite{Roos2008,Stynes2005}.

Physics-informed neural networks provide a residual-based framework for approximating solutions of differential equations. In the standard formulation, the unknown solution is represented by a neural network, and the trainable parameters are obtained by minimizing a loss function that contains the differential-equation residual together with penalty terms for the boundary or initial conditions \cite{Raissi2019,Karniadakis2021}. Derivatives of the network output are computed by automatic differentiation, which makes the formulation flexible for different differential operators and boundary-value problems \cite{Lu2021}.

Despite this flexibility, the standard formulation can face serious training and approximation difficulties. The optimization problem is nonlinear and non-convex with respect to the network parameters, so the final approximation can depend on initialization, architecture, optimizer choice, collocation distribution, and loss scaling. Moreover, in the usual soft-boundary setting, the boundary conditions are imposed as penalty terms in the loss function rather than being satisfied analytically. This creates competing objectives: the network must reduce the residual inside the domain while also learning the boundary constraints. Gradient pathologies, failure modes, loss-landscape issues, and training difficulties of such formulations have been studied in several works \cite{Wang2021,Krishnapriyan2021,Basir2022,Wang2022NTK}. In particular, De Florio et al. \cite{DeFlorio2024} reported that standard residual-based formulations can become unreliable for second-order ordinary differential equations with sharp gradients.

Several approaches have been proposed to reduce some limitations of standard residual-based PINN formulations. Constrained expressions and the Theory of Functional Connections can be used to satisfy prescribed constraints analytically \cite{Leake2020,Mortari2017}. Extreme-learning-machine-based approaches reduce the training burden by fixing hidden-layer parameters and solving mainly for output weights \cite{Dwivedi2020,Huang2006}. X-TFC combines constrained expressions with extreme-learning ideas \cite{Schiassi2021}. Domain-decomposed and variational formulations have also been developed to improve the treatment of complex differential problems \cite{Jagtap2020,Kharazmi2021}. These developments show that the standard soft-boundary PINN used in this work is a baseline formulation, not a representative of all physics-informed neural methods.

The present paper takes a different route. Instead of proposing another neural variant, it uses Legendre wavelet collocation as a local algebraic benchmark against the standard soft-boundary PINN. This comparison is useful because the two methods use different approximation mechanisms. The standard PINN seeks the solution through residual-loss minimization over trainable network parameters, whereas the Legendre wavelet method builds a local polynomial approximation and determines its coefficients from collocation equations.

Legendre wavelet methods are attractive in this setting because their basis functions are locally supported on sub-intervals. A global polynomial basis acts over the whole domain, whereas a local wavelet basis acts only on its own cell. This locality is useful for boundary-layer problems, where the important variation occurs only near a small part of the interval. Legendre wavelet direct and operational-matrix methods have been used to transform differential equations into algebraic systems in several settings \cite{Razzaghi2000,Jafari2011,Balaji2014,Rehman2011,Nemati2019}.

The gap addressed here is the absence of a direct benchmark between a standard soft-boundary PINN and a local Legendre wavelet collocation method on the same singularly perturbed transport problem. Existing studies mainly examine residual-trained neural formulations or modifications of the standard PINN framework, whereas the present work places a baseline residual-loss PINN and a cellwise algebraic wavelet method under the same exact-solution test. This comparison isolates two different approximation mechanisms: global residual-loss minimization through trainable network parameters and local algebraic collocation through Legendre wavelet coefficients. The main contribution is the construction and testing of a direct two-sum Legendre wavelet collocation scheme against the standard soft-boundary PINN for $\mathrm{Pe}=1,10,100,$ and $1000$. The comparison uses dense-grid error evaluation with additional points near the boundary layer, so that errors hidden by coarse sampling can be detected.

The rest of the paper is organized as follows. Section 2 describes the benchmark problem and its exact solution. Section 3 presents the direct two-sum Legendre wavelet collocation method. Section 4 describes the standard soft-boundary PINN formulation. Section 5 reports the numerical results and comparison. Section 6 concludes the paper.

\section{Benchmark Problem}

We consider the steady one-dimensional advection--diffusion boundary-value problem
\begin{equation}
	u''(x)-\mathrm{Pe} \, u'(x)=0,\qquad 0\leq x\leq 1,
	\label{eq:bvp}
\end{equation}
subject to the boundary conditions
\begin{equation}
	u(0)=1,\qquad u(1)=0.
	\label{eq:bc}
\end{equation}
Here $\mathrm{Pe}$ denotes the Peclet number and is related to the diffusion coefficient $\nu$ by
\begin{equation}
	\mathrm{Pe}=\frac{1}{\nu}.
	\label{eq:peclet}
\end{equation}
In this non-dimensional form, decreasing $\nu$ increases $\mathrm{Pe}$ and produces a thin boundary layer near the right boundary.

The exact solution of Eqs.~\eqref{eq:bvp}--\eqref{eq:bc} is
\begin{equation}
	u_{\mathrm{exact}}(x)=
	\frac{1-\exp\left(\mathrm{Pe}(x-1)\right)}
	{1-\exp(-\mathrm{Pe})}.
	\label{eq:exact}
\end{equation}
For small $\mathrm{Pe}$, the solution varies smoothly across the interval. For large $\mathrm{Pe}$, the solution remains close to one over most of the domain and then rapidly decreases near $x=1$ in order to satisfy the boundary condition $u(1)=0$. Therefore, the main numerical difficulty is the accurate resolution of the narrow boundary layer near the right endpoint.

In this work, four test cases are considered:
\begin{equation}
	\nu=1,\;0.1,\;0.01,\;0.001,
\end{equation}
or equivalently,
\begin{equation}
	\mathrm{Pe}=1,\;10,\;100,\;1000.
\end{equation}
These cases represent a progression from a smooth solution to a strongly boundary-layer-dominated solution. The same exact solution, dense evaluation grid, and error metrics are used for both numerical methods so that the comparison remains consistent.

\section{Legendre Wavelet Collocation Method}
This section describes the direct two-sum Legendre wavelet collocation method used in this work. The term ``direct'' indicates that the approximate solution itself is expanded in Legendre wavelet basis functions and the differential equation is enforced at collocation points, without first converting the problem into an integral equation. The term ``two-sum'' refers to the two indices in the approximation: one for the subinterval and one for the Legendre mode inside that subinterval. Since the basis functions are locally supported, the approximation is constructed cell by cell, which is suitable for the boundary-layer structure considered here. Legendre wavelet methods have been used in several differential-equation problems because they combine local polynomial approximation with the reduction of differential problems to algebraic systems \cite{Razzaghi2000,Jafari2011,Balaji2014,Rehman2011,Nemati2019,raygupta2018book,NayakGupta2026Heat,NayakGupta2025}.

For the boundary-value problem in Eq.~\eqref{eq:bvp}, we divide the interval $[0,1]$ into $J$ equal cells. Let
\begin{equation}
	x_j=\frac{j}{J},\qquad j=0,1,\ldots,J,
\end{equation}
and define
\begin{equation}
	I_j=[x_j,x_{j+1}],\qquad j=0,1,\ldots,J-1.
\end{equation}
The length of each cell is
\begin{equation}
	h=\frac{1}{J}.
\end{equation}
For a point $x\in I_j$, the physical cell is mapped to the standard Legendre interval $[-1,1]$ by
\begin{equation}
	z_j(x)=\frac{2(x-x_j)}{h}-1.
\end{equation}
Thus, each cell is treated using the same local coordinate $z$.

Let $L_m(z)$ denote the Legendre polynomial of degree $m$ on $[-1,1]$. The local Legendre wavelet basis function on cell $I_j$ is defined by

\begin{equation}
	\psi_{j,m}(x)=
	\begin{cases}
		\sqrt{\dfrac{2m+1}{h}}\,L_m\!\left(z_j(x)\right), & x\in I_j,\\[2mm]
		0, & x\notin I_j,
	\end{cases}
	\label{eq:lwm_basis}
\end{equation}

where $m=0,1,\ldots,M-1$. Here $M$ is the number of Legendre modes used in each cell. The factor $\sqrt{(2m+1)/h}$ is the standard normalization factor for the mapped Legendre basis. The important point is that $\psi_{j,m}$ is a local basis function: it is zero outside its own cell. At internal interfaces, the values from the left and right cells are treated as one-sided traces.

The approximate solution of Eq.~\eqref{eq:bvp} is written in the finite two-sum form as
\begin{equation}
	u_{J,M}(x)=
	\sum_{j=0}^{J-1}\sum_{m=0}^{M-1}
	c_{j,m}\psi_{j,m}(x),
	\label{eq:lwm_approx}
\end{equation}
where $c_{j,m}$ are unknown coefficients. The first sum runs over the cells, and the second sum runs over the Legendre modes inside each cell. Therefore, the total number of unknown coefficients is
\begin{equation}
	N_{\mathrm{unk}}=JM.
\end{equation}

Since the governing equation contains $u'(x)$ and $u''(x)$, the derivatives of the approximation are obtained by differentiating the local basis functions. From the mapping between $x$ and $z$,

\begin{equation}
	\frac{dz_j}{dx}=\frac{2}{h},
\end{equation}
\begin{equation}
	\psi'_{j,m}(x)=
	\sqrt{\frac{2m+1}{h}}\frac{2}{h}
	L'_m\!\left(z_j(x)\right),
	\qquad x\in I_j,
\end{equation}

and the second derivative 
\begin{equation}
	\psi''_{j,m}(x)=
	\sqrt{\frac{2m+1}{h}}
	\left(\frac{2}{h}\right)^2
	L''_m\!\left(z_j(x)\right),
	\qquad x\in I_j.
\end{equation}
Consequently,
\begin{equation}
	u'_{J,M}(x)=
	\sum_{j=0}^{J-1}\sum_{m=0}^{M-1}
	c_{j,m}\psi'_{j,m}(x),
\end{equation}
and
\begin{equation}
	u''_{J,M}(x)=
	\sum_{j=0}^{J-1}\sum_{m=0}^{M-1}
	c_{j,m}\psi''_{j,m}(x).
\end{equation}

For the benchmark problem, the residual is defined by

\begin{equation}
	R(x)=u''_{J,M}(x)-\mathrm{Pe}\,u'_{J,M}(x).
	\label{eq:lwm_residual}
\end{equation}

The collocation idea is to force this residual to vanish at selected points inside each cell. To obtain a square algebraic system, the residual equations are combined with the two boundary equations and the interface-continuity equations. Therefore, we use $M-2$ interior collocation points in each cell. These points are chosen as Gauss--Legendre points on $[-1,1]$ and then mapped to the physical cell. If $\xi_r$, $r=1,2,\ldots,M-2$, are the Gauss--Legendre points, then the corresponding physical collocation points in cell $I_j$ are

\begin{equation}
	x_{j,r}=x_j+\frac{h}{2}(1+\xi_r),
	\qquad r=1,2,\ldots,M-2.
\end{equation}
At these points, the residual equations are
\begin{equation}
	R(x_{j,r})=0,
	\qquad j=0,1,\ldots,J-1,\quad r=1,2,\ldots,M-2.
	\label{eq:lwm_residual_equations}
\end{equation}
This gives $J(M-2)$ residual equations.

The two boundary conditions are imposed directly:

\begin{equation}
	u_{J,M}(0)=1,
	\qquad
	u_{J,M}(1)=0.
	\label{eq:lwm_bc_equations}
\end{equation}
These two equations enforce the prescribed endpoint values.

Since the approximation is constructed cell by cell, neighbouring cells must be connected at their common interfaces. The internal interface points are
\begin{equation}
	x_s=\frac{s}{J},
	\qquad s=1,2,\ldots,J-1.
\end{equation}
At each interface, we impose continuity of the solution,
\begin{equation}
	u_{J,M}(x_s^-)=u_{J,M}(x_s^+),
	\label{eq:lwm_value_continuity}
\end{equation}
and continuity of the first derivative,
\begin{equation}
	u'_{J,M}(x_s^-)=u'_{J,M}(x_s^+).
	\label{eq:lwm_derivative_continuity}
\end{equation}

The first condition prevents jumps in the numerical solution. The second condition prevents jumps in the slope, which is natural because the differential equation contains a second derivative. Since there are $J-1$ internal interfaces, Eqs.~\eqref{eq:lwm_value_continuity}--\eqref{eq:lwm_derivative_continuity} provide $2(J-1)$ equations.

To assemble the algebraic system, the coefficients are arranged in the vector

\begin{equation}
	\mathbf{c}=
	\left(c_{0,0},c_{0,1},\ldots,c_{0,M-1},
	c_{1,0},\ldots,c_{J-1,M-1}\right)^T.
\end{equation}

The column corresponding to the coefficient $c_{j,m}$ is denoted by
\begin{equation}
	p(j,m)=jM+m.
\end{equation}

Let $\ell$ denote the row index of the algebraic system.

For a residual collocation point $x_{j,r}$, the corresponding row is

\begin{equation}
	A_{\ell,p(k,m)}
	=
	\psi''_{k,m}(x_{j,r})
	-
	\mathrm{Pe}\,\psi'_{k,m}(x_{j,r}),
	\qquad
	k=0,1,\ldots,J-1,\quad m=0,1,\ldots,M-1.
	\label{eq:lwm_residual_matrix_row}
\end{equation}
with
\begin{equation}
	b_{\ell}=0.
\end{equation}
Because the basis functions are locally supported, only the basis functions belonging to the active cell contribute at an interior collocation point.

The left boundary condition gives
\begin{equation}
	A_{\ell,p(0,m)}=\psi_{0,m}(0),
	\qquad m=0,1,\ldots,M-1,
\end{equation}
with all other entries in that row equal to zero and
\begin{equation}
	b_{\ell}=1.
\end{equation}
Similarly, the right boundary condition gives
\begin{equation}
	A_{\ell,p(J-1,m)}=\psi_{J-1,m}(1),
	\qquad m=0,1,\ldots,M-1,
\end{equation}
with all other entries in that row equal to zero and
\begin{equation}
	b_{\ell}=0.
\end{equation}

At an internal interface $x_s=s/J$, the value-continuity row is
\begin{equation}
	A_{\ell,p(s-1,m)}=\psi_{s-1,m}(x_s),
	\qquad
	A_{\ell,p(s,m)}=-\psi_{s,m}(x_s),
	\qquad m=0,1,\ldots,M-1,
\end{equation}
with
\begin{equation}
	b_{\ell}=0.
\end{equation}
The derivative-continuity row is
\begin{equation}
	A_{\ell,p(s-1,m)}=\psi'_{s-1,m}(x_s),
	\qquad
	A_{\ell,p(s,m)}=-\psi'_{s,m}(x_s),
	\qquad m=0,1,\ldots,M-1,
\end{equation}
again with
\begin{equation}
	b_{\ell}=0.
\end{equation}

Thus, the residual, boundary, and interface-continuity equations are assembled into one algebraic system for the coefficient vector $\mathbf{c}$. The total number of equations is
\begin{equation}
	J(M-2)+2+2(J-1)=JM.
\end{equation}
This is exactly equal to the number of unknown coefficients in Eq.~\eqref{eq:lwm_approx}. Therefore, the collocation procedure produces the square linear algebraic system
\begin{equation}
	A\mathbf{c}=\mathbf{b}.
	\label{eq:lwm_linear_system}
\end{equation}
After solving this system, the numerical solution is reconstructed from Eq.~\eqref{eq:lwm_approx}.

The algebraic residual
\begin{equation}
	\left\|A\mathbf{c}-\mathbf{b}\right\|_{\infty}
\end{equation}
is reported to check the linear solve, while the pointwise error is computed against the exact solution with additional grid resolution near $x=1$.

\section{Standard Soft-Boundary PINN Formulation}
In this section, we describe the standard soft-boundary PINN formulation used for comparison with the direct Legendre wavelet collocation method. The boundary conditions are not built into the trial solution. Instead, the differential equation and boundary conditions are imposed through the differential-equation residual and boundary penalty terms in the loss function. No hard-boundary construction, adaptive loss weighting, domain decomposition, or boundary-layer correction is used, so that the comparison reflects the behavior of the baseline PINN formulation.

Let $u_{\theta}(x)$ denote the neural approximation of the unknown solution $u(x)$, where $\theta$ represents all trainable weights and biases. The network takes the spatial coordinate $x\in[0,1]$ as input and returns one scalar output $u_{\theta}(x)$. In the numerical experiments, a fully connected feed-forward neural network is used with one input neuron, four hidden layers, 64 neurons in each hidden layer, the hyperbolic tangent activation function, and one output neuron. After training, $u_{\theta}(x)$ is used as the PINN approximation of the solution.

The benchmark equation contains both first- and second-order derivatives. Therefore, the derivatives of $u_{\theta}(x)$ with respect to $x$ are required. These derivatives are computed using automatic differentiation:
\begin{equation}
	u'_{\theta}(x)=\frac{d u_{\theta}(x)}{dx},
	\qquad
	u''_{\theta}(x)=\frac{d^2 u_{\theta}(x)}{dx^2}.
\end{equation}

Using these derivatives, the differential-equation residual is defined as
\begin{equation}
	r_{\theta}(x)=u''_{\theta}(x)-\mathrm{Pe}\,u'_{\theta}(x).
	\label{eq:pinn_residual}
\end{equation}
If $r_{\theta}(x)=0$ in the interval and the boundary conditions are satisfied, then $u_{\theta}(x)$ satisfies the boundary-value problem.

The residual is enforced at a finite set of collocation points. Let
\begin{equation}
	\left\{x_i^r\right\}_{i=1}^{N_r}\subset[0,1],
\end{equation}
denote the residual collocation points. In this work, uniformly distributed collocation points are used with $N_r=10000$. The residual loss is defined by
\begin{equation}
	\mathcal{L}_r(\theta)=
	\frac{1}{N_r}
	\sum_{i=1}^{N_r}
	\left|r_{\theta}(x_i^r)\right|^2 .
	\label{eq:pinn_residual_loss}
\end{equation}
This term measures how well the neural approximation satisfies the differential equation on the chosen collocation set.

The boundary conditions are imposed through the boundary loss
\begin{equation}
	\mathcal{L}_b(\theta)=
	\left|u_{\theta}(0)-1\right|^2+
	\left|u_{\theta}(1)\right|^2 .
	\label{eq:pinn_boundary_loss}
\end{equation}
The first term enforces $u(0)=1$, while the second term enforces $u(1)=0$.

The total loss function is taken as
\begin{equation}
	\mathcal{L}(\theta)=
	\mathcal{L}_r(\theta)+\lambda_b\mathcal{L}_b(\theta),
	\label{eq:pinn_total_loss}
\end{equation}
where $\lambda_b$ is the boundary-loss weight. In the computations, we use
\begin{equation}
	\lambda_b=100.
\end{equation}
The trainable parameters are obtained by minimizing the total loss:
\begin{equation}
	\theta^{\ast}=\arg\min_{\theta}\mathcal{L}(\theta).
\end{equation}

The optimization is performed in two stages. First, the Adam optimizer is used for $15000$ epochs with learning rate $10^{-3}$. Then, the L-BFGS optimizer is applied as a refinement step with a maximum of $500$ iterations. This two-stage strategy is used because Adam provides a stable initial descent, while L-BFGS can further refine the residual minimization.

This is a soft-boundary formulation because the boundary conditions are imposed as penalty terms in the loss function. Therefore, the approximation is not guaranteed to satisfy $u(0)=1$ and $u(1)=0$ exactly after training. This point is important for sharp-gradient problems. When the boundary layer becomes very thin, the optimization problem must balance the interior residual loss and the boundary loss. A small residual loss alone does not necessarily imply that the boundary-layer profile has been accurately captured.

After training, the numerical solution is defined by
\begin{equation}
	u_{\mathrm{PINN}}(x)=u_{\theta^*}(x).
\end{equation}
The trained approximation is compared with the exact solution on a dense evaluation grid. This grid is independent of the training collocation points and includes additional points near $x=1$, where the boundary layer occurs. The reported error measures are the maximum absolute error, mean absolute error, root-mean-square error, an $L^2$-type error, and the absolute boundary errors at $x=0$ and $x=1$.

\section{Numerical Results and Discussion}

\subsection{Computational setup and error measures}

The numerical experiments are performed for four Péclet numbers,
\begin{equation}
	\mathrm{Pe}=1,\;10,\;100,\;1000,
\end{equation}
corresponding to
\begin{equation}
	\nu=1,\;0.1,\;0.01,\;0.001,
\end{equation}
since $\mathrm{Pe}=1/\nu$. These cases represent a gradual transition from a smooth solution to a strongly boundary-layer-dominated solution. The exact solution in Eq.~\eqref{eq:exact} is used as the reference solution in all comparisons.

For a numerical approximation $u_{\mathrm{num}}(x)$, the pointwise absolute error at an evaluation point $x_i$ is defined by
\begin{equation}
	e(x_i)=\left|u_{\mathrm{num}}(x_i)-u_{\mathrm{exact}}(x_i)\right|.
\end{equation}
The maximum absolute error is
\begin{equation}
	E_{\infty}=\max_{1\leq i\leq N} e(x_i),
\end{equation}
the grid-based mean absolute error is
\begin{equation}
	E_{\mathrm{mean}}=\frac{1}{N}\sum_{i=1}^{N}e(x_i),
\end{equation}
and the grid-based root-mean-square error is
\begin{equation}
	E_{\mathrm{rms}}=
	\left(\frac{1}{N}\sum_{i=1}^{N}e(x_i)^2\right)^{1/2}.
\end{equation}
An $L^2$-type error is also reported and is approximated numerically on the evaluation grid:
\begin{equation}
	E_2=
	\left(\int_0^1
	\left|u_{\mathrm{num}}(x)-u_{\mathrm{exact}}(x)\right|^2\,dx
	\right)^{1/2}.
\end{equation}
All error measures are evaluated on a dense grid that combines uniformly distributed points with additional points clustered near $x=1$, where the boundary layer occurs. This is important because a coarse uniform grid can miss the main pointwise error for large $\mathrm{Pe}$.
\subsection{Direct two-sum Legendre wavelet results}

The direct two-sum Legendre wavelet method is tested using $M=6$ Legendre modes in each cell. The number of cells is increased with the Peclet number in order to resolve the increasingly sharp boundary layer. The values used are $J=4,8,40,$ and $200$ for $\mathrm{Pe}=1,10,100,$ and $1000$, respectively.

\begin{figure}[H]
	\centering
	\begin{minipage}{0.48\textwidth}
		\centering
		\includegraphics[width=\textwidth]{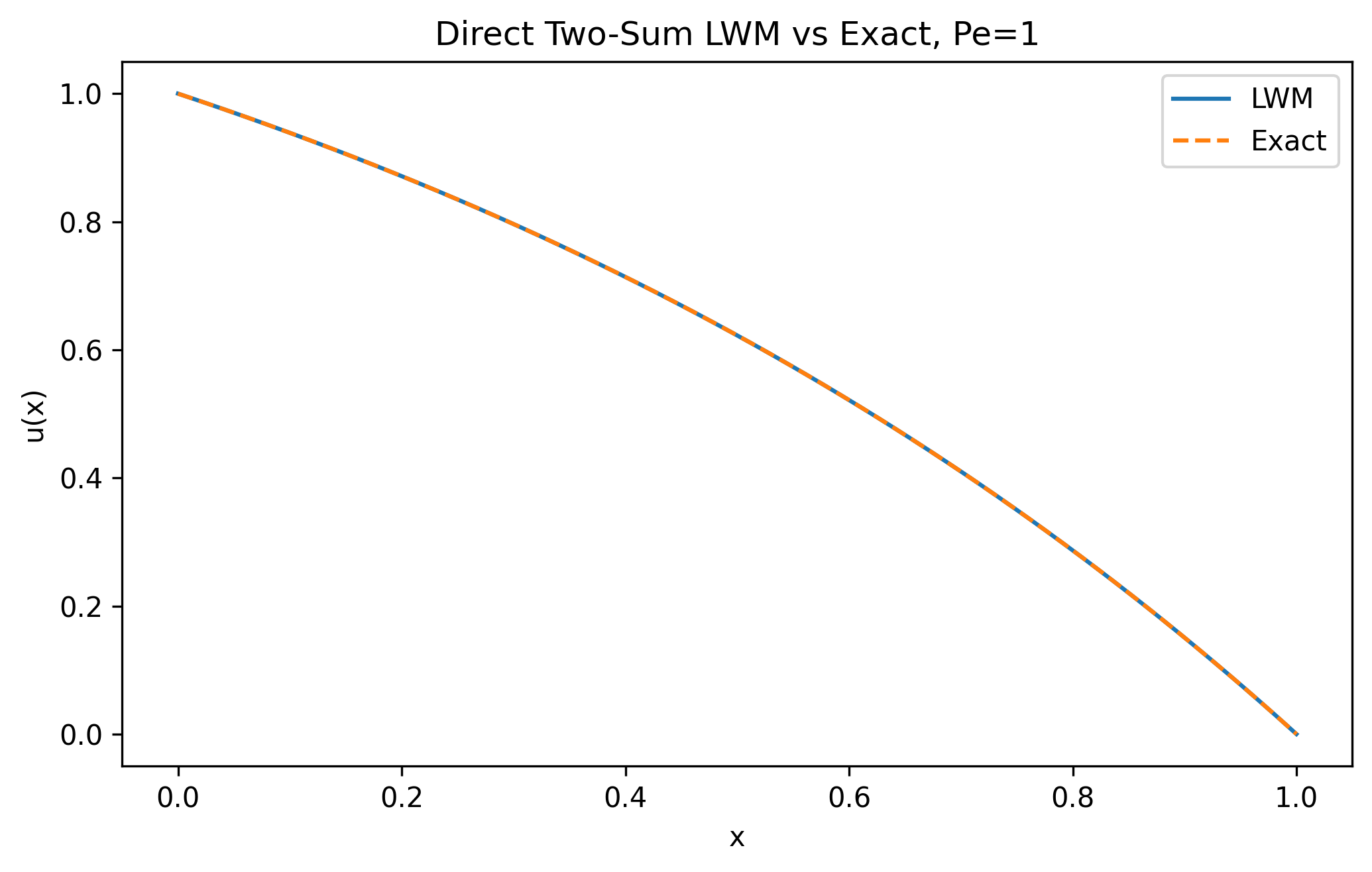}\\
		\textbf{(a)} $\mathrm{Pe}=1$
	\end{minipage}
	\begin{minipage}{0.48\textwidth}
		\centering
		\includegraphics[width=\textwidth]{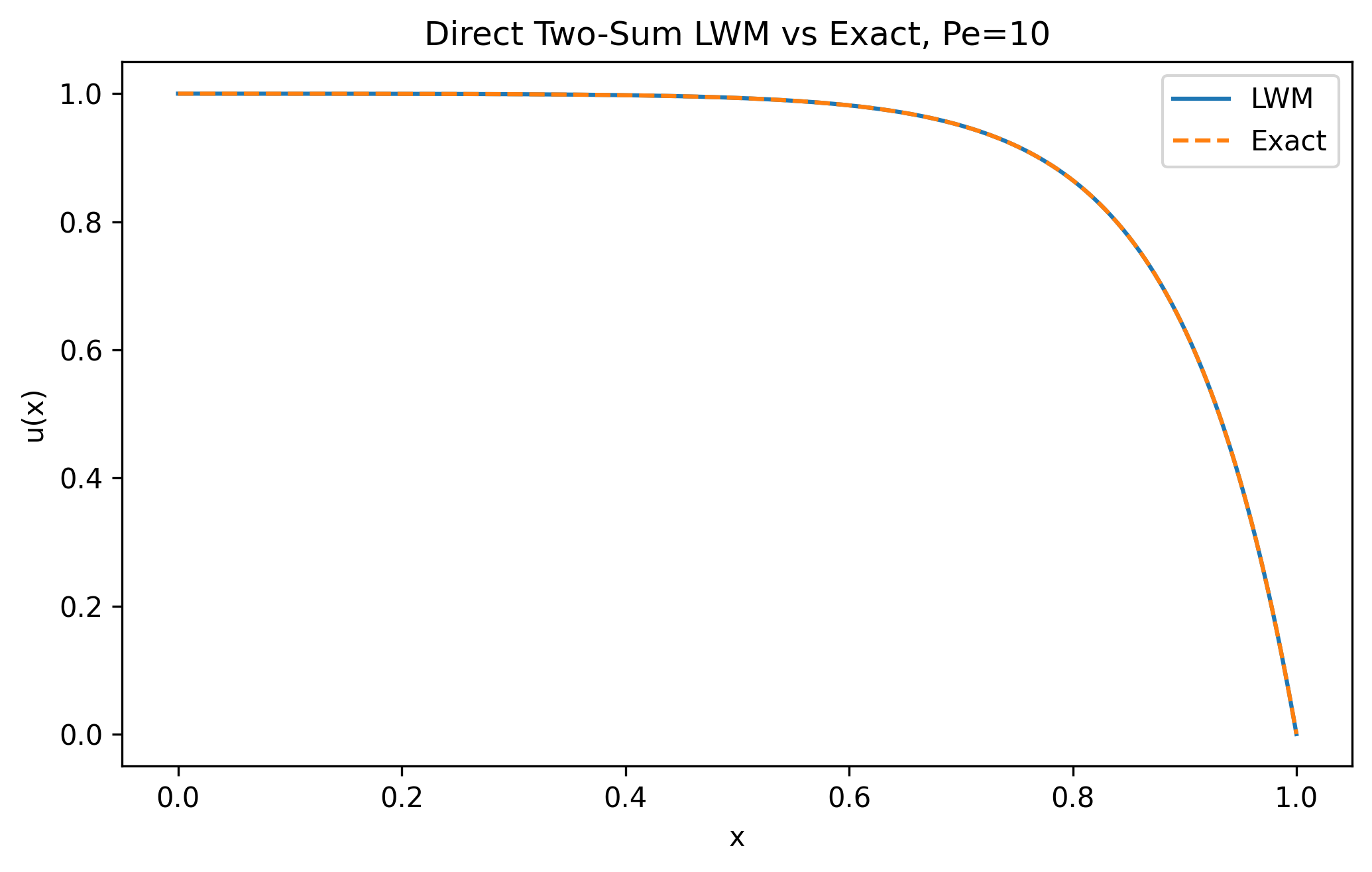}\\
		\textbf{(b)} $\mathrm{Pe}=10$
	\end{minipage}
	
	\vspace{2mm}
	
	\begin{minipage}{0.48\textwidth}
		\centering
		\includegraphics[width=\textwidth]{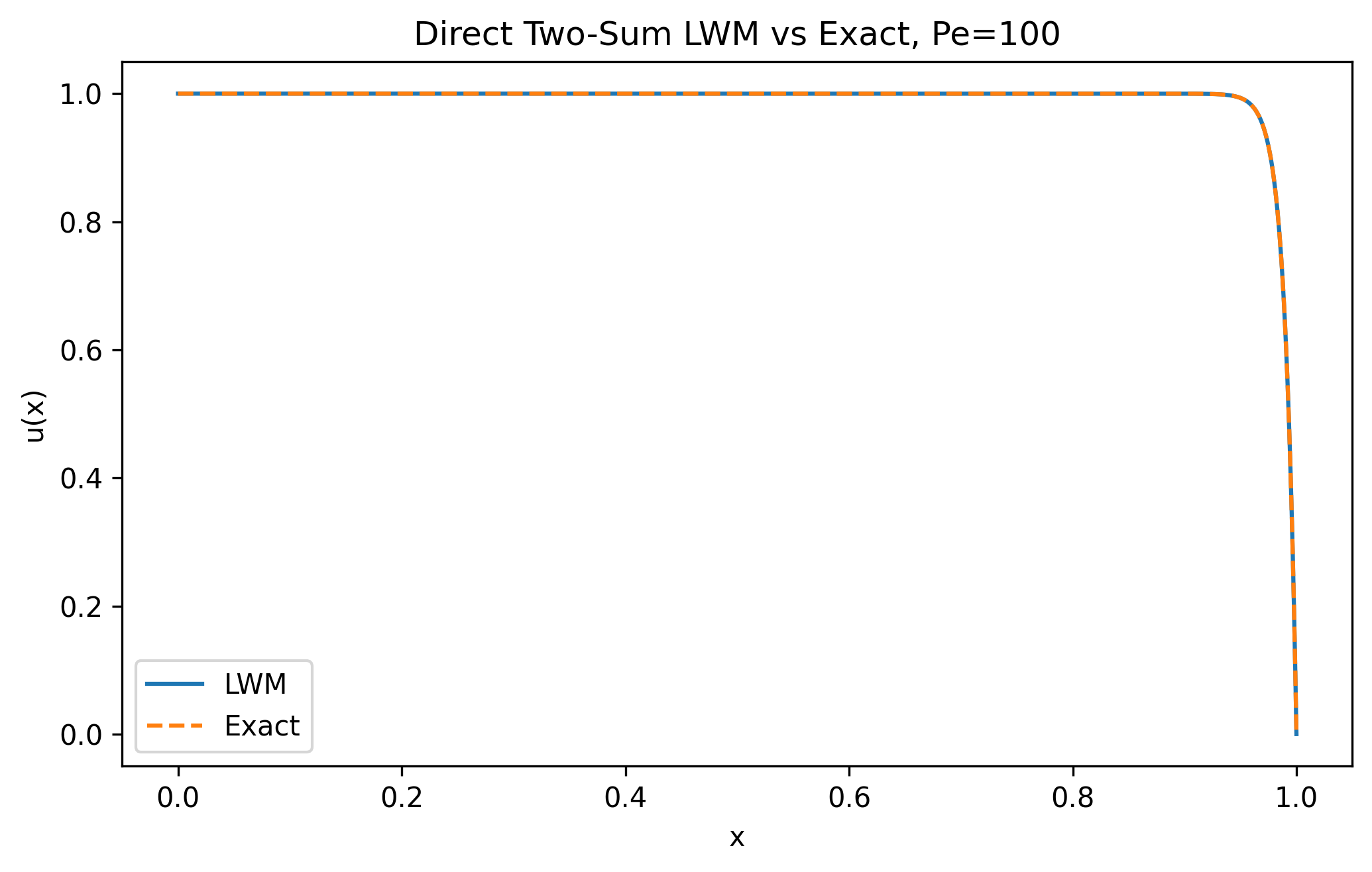}\\
		\textbf{(c)} $\mathrm{Pe}=100$
	\end{minipage}
	\begin{minipage}{0.48\textwidth}
		\centering
		\includegraphics[width=\textwidth]{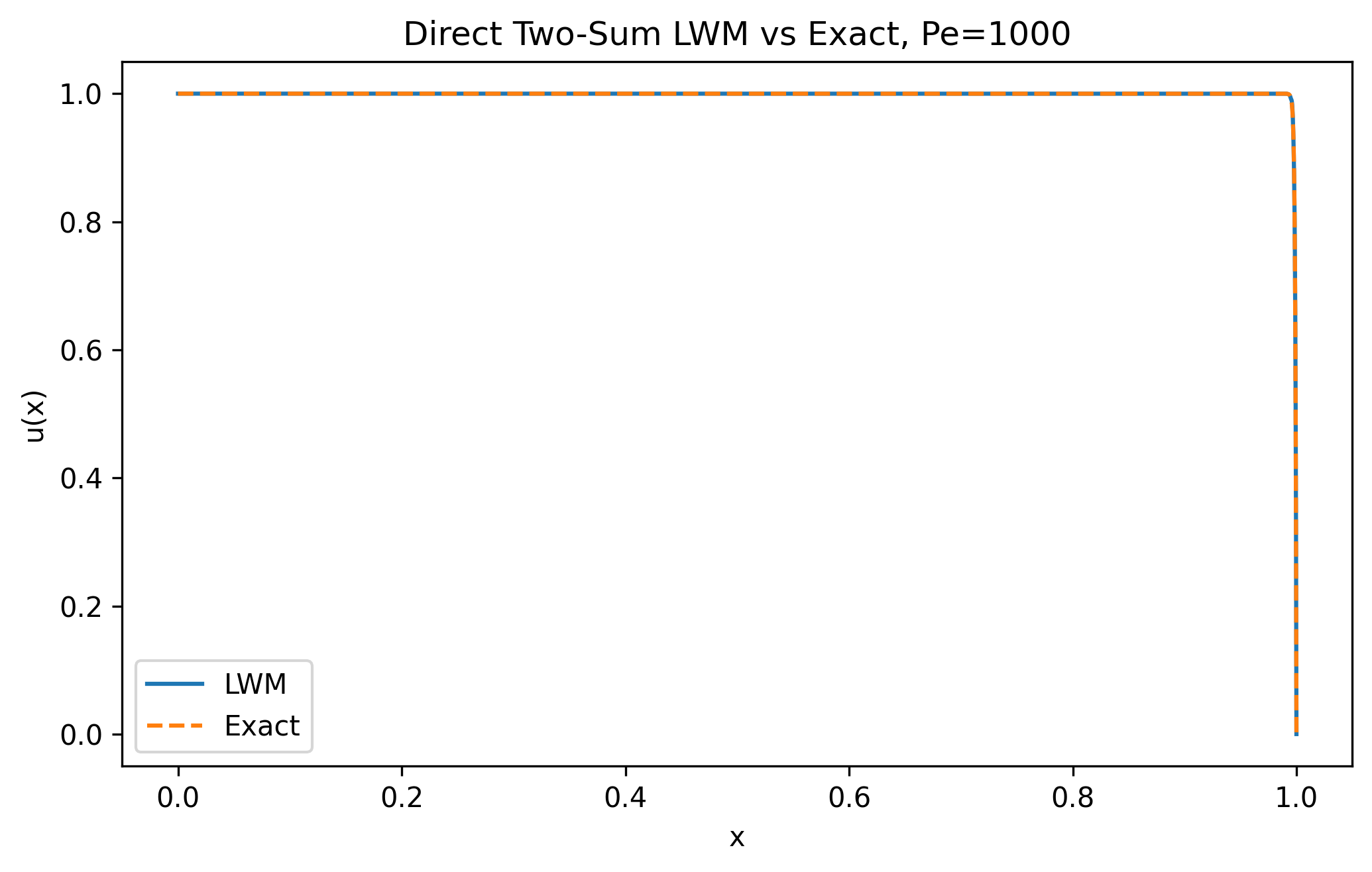}\\
		\textbf{(d)} $\mathrm{Pe}=1000$
	\end{minipage}
	\caption{Direct two-sum Legendre wavelet solution and exact solution for different Peclet numbers.}
	\label{fig:lwm_solution_all}
\end{figure}

Figure~\ref{fig:lwm_solution_all} shows the LWM solution together with the exact solution for all four Peclet numbers. At the plotted scale, the LWM and exact curves are visually very close for all cases. For larger $\mathrm{Pe}$, the solution remains close to one over most of the domain and then changes rapidly near $x=1$.

\begin{figure}[H]
	\centering
	\begin{minipage}{0.48\textwidth}
		\centering
		\includegraphics[width=\textwidth]{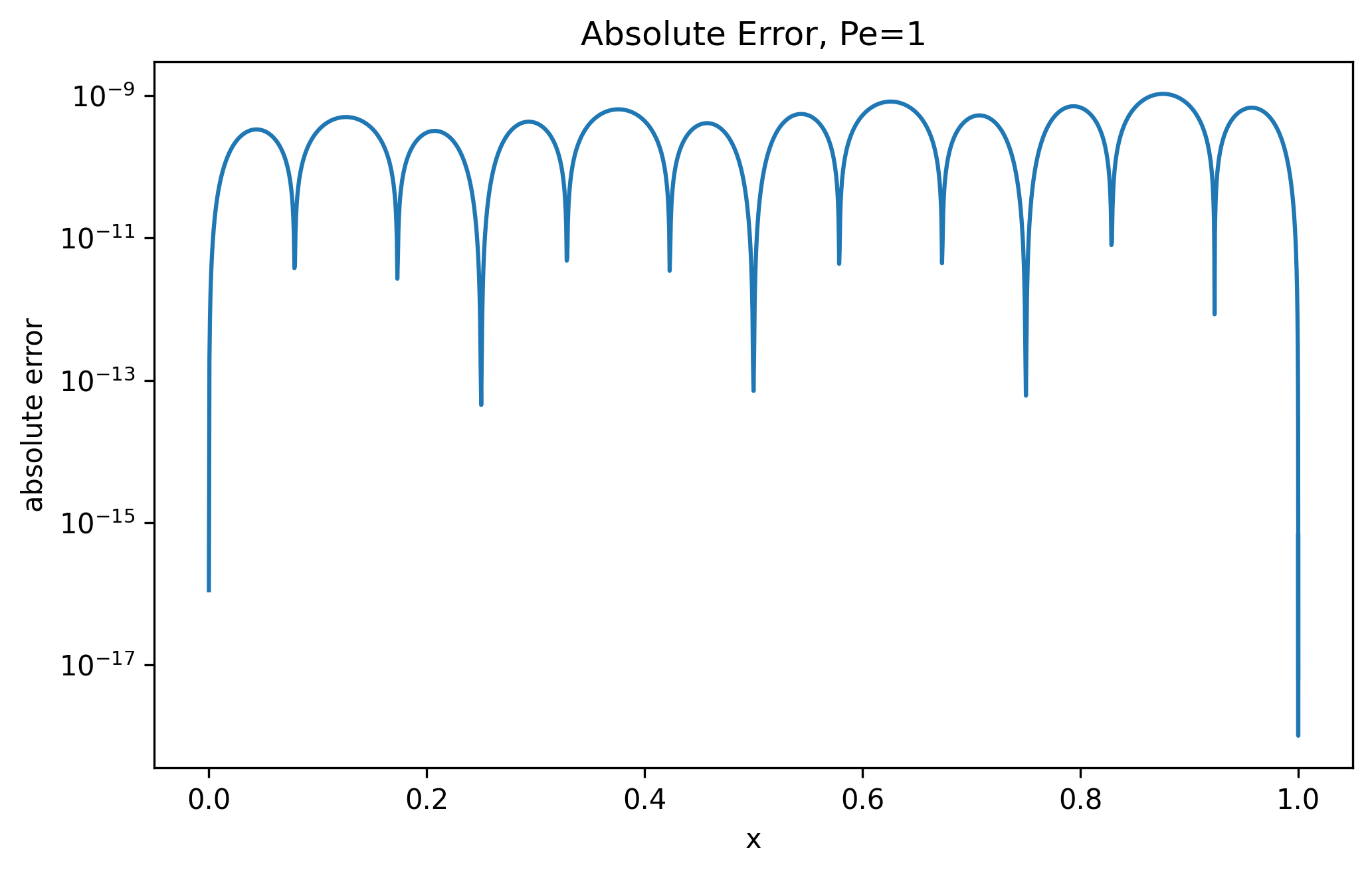}\\
		\textbf{(a)} $\mathrm{Pe}=1$
	\end{minipage}
	\begin{minipage}{0.48\textwidth}
		\centering
		\includegraphics[width=\textwidth]{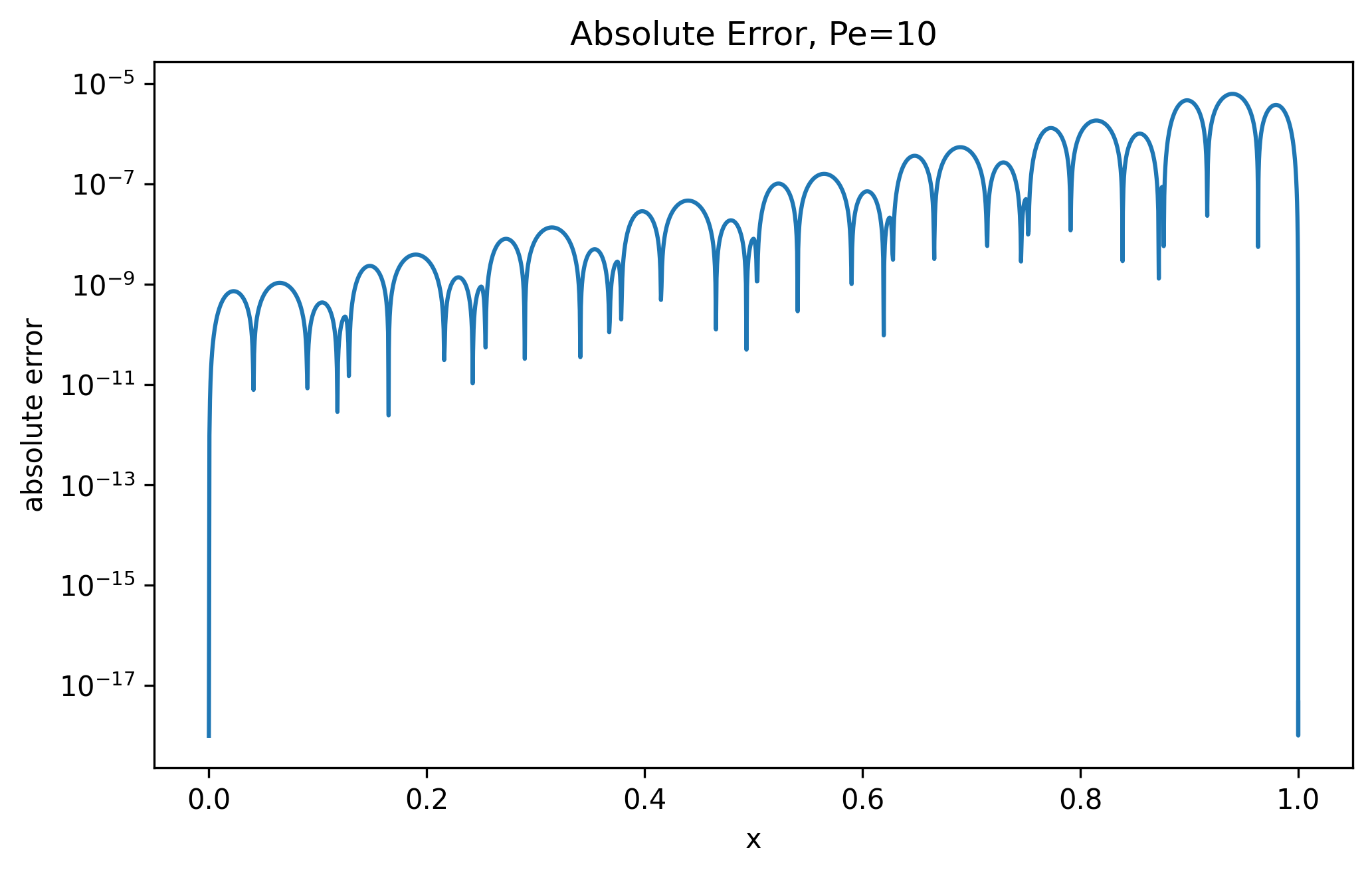}\\
		\textbf{(b)} $\mathrm{Pe}=10$
	\end{minipage}
	
	\vspace{2mm}
	
	\begin{minipage}{0.48\textwidth}
		\centering
		\includegraphics[width=\textwidth]{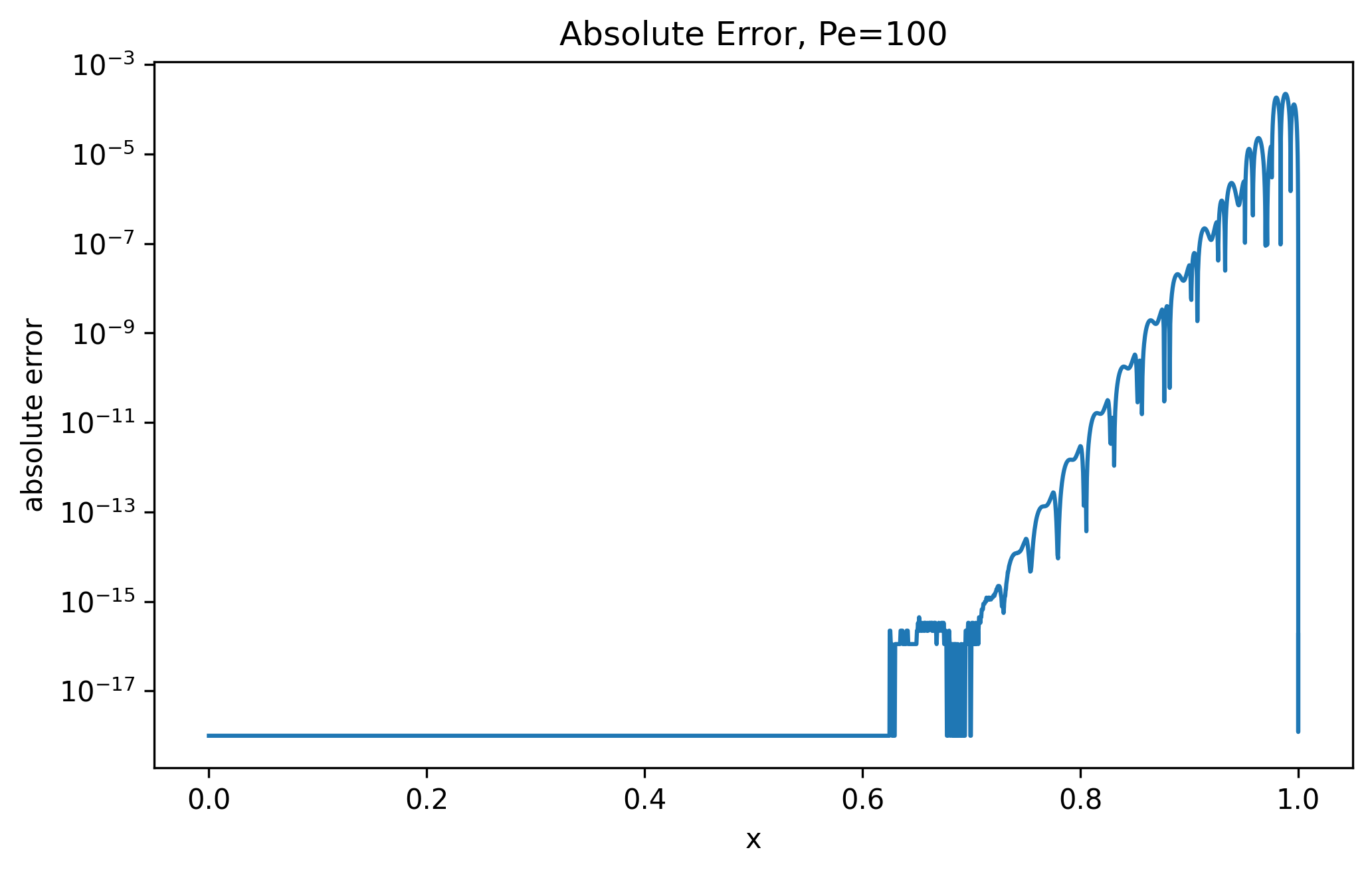}\\
		\textbf{(c)} $\mathrm{Pe}=100$
	\end{minipage}
	\begin{minipage}{0.48\textwidth}
		\centering
		\includegraphics[width=\textwidth]{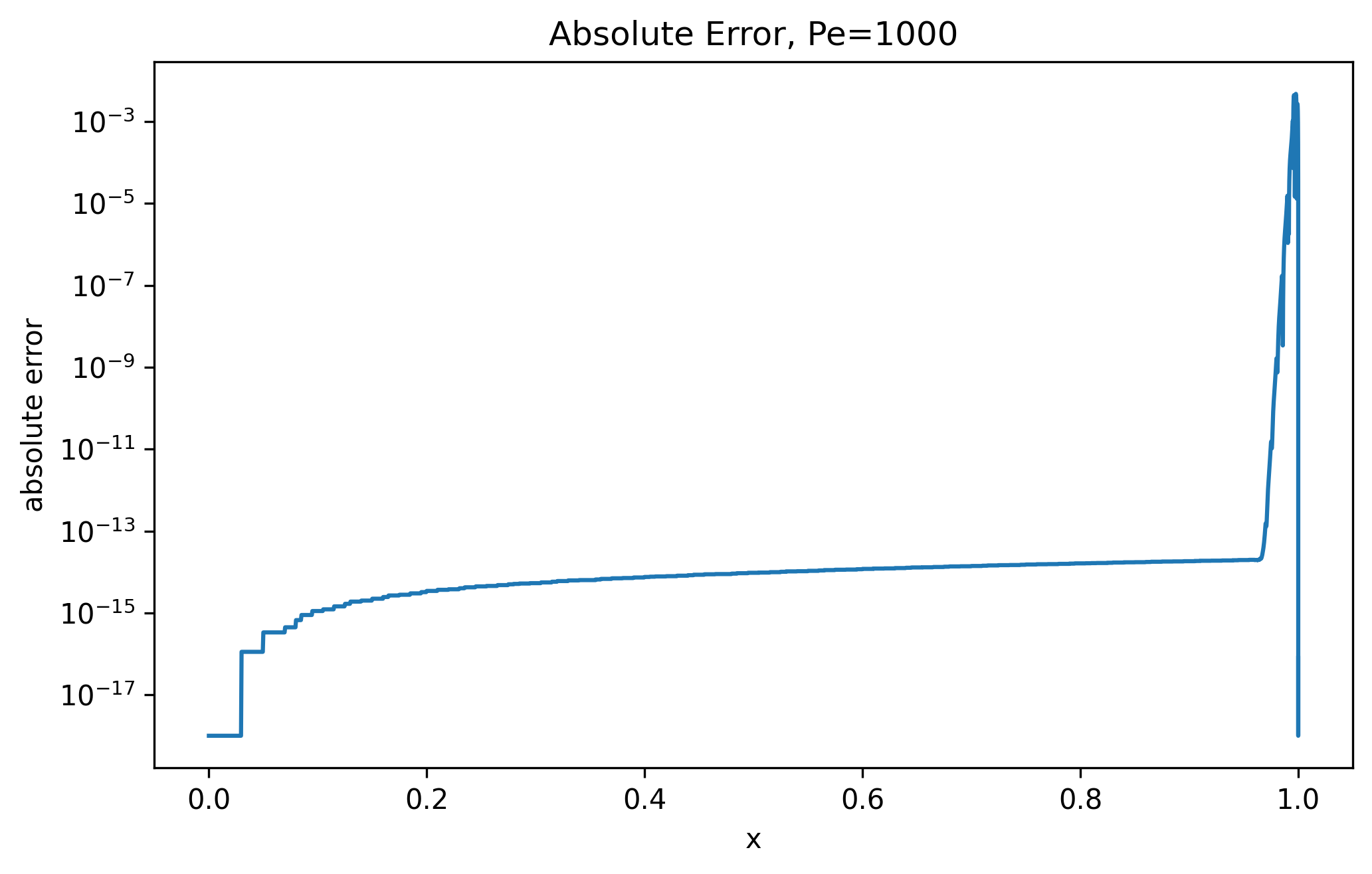}\\
		\textbf{(d)} $\mathrm{Pe}=1000$
	\end{minipage}
	\caption{Absolute error of the direct two-sum Legendre wavelet method for different Peclet numbers.}
	\label{fig:lwm_error_all}
\end{figure}

Figure~\ref{fig:lwm_error_all} shows the corresponding absolute-error profiles. The errors are small for $\mathrm{Pe}=1$ and $\mathrm{Pe}=10$. For $\mathrm{Pe}=100$ and $\mathrm{Pe}=1000$, the largest errors occur near the sharp transition close to $x=1$. This is consistent with the increasing boundary-layer difficulty as $\mathrm{Pe}$ increases.

The pointwise LWM results are reported in Tables~\ref{tab:lwm_pe1_pointwise}--\ref{tab:lwm_pe1000_pointwise}. The tabulated points include standard points in the interval and additional points near $x=1$ for the sharper boundary-layer cases.

\begin{table}[H]
	\centering
	\scriptsize
	\caption{Pointwise LWM results for $\mathrm{Pe}=1$.}
	\label{tab:lwm_pe1_pointwise}
	\begin{tabular}{||lllllll||}
		\hline\hline
		$x$ & & LWM solution & & Exact solution & & Absolute error \\
		\hline\hline
		0.0000E+00 & & 1.0000E+00 & & 1.0000E+00 & & 1.1102E-16 \\
		1.0000E-01 & & 9.3879E-01 & & 9.3879E-01 & & 3.1972E-10 \\
		2.0000E-01 & & 8.7115E-01 & & 8.7115E-01 & & 3.0143E-10 \\
		3.0000E-01 & & 7.9639E-01 & & 7.9639E-01 & & 4.1222E-10 \\
		4.0000E-01 & & 7.1377E-01 & & 7.1377E-01 & & 4.3892E-10 \\
		5.0000E-01 & & 6.2246E-01 & & 6.2246E-01 & & 7.0721E-14 \\
		6.0000E-01 & & 5.2155E-01 & & 5.2155E-01 & & 5.2718E-10 \\
		7.0000E-01 & & 4.1002E-01 & & 4.1002E-01 & & 4.9697E-10 \\
		8.0000E-01 & & 2.8676E-01 & & 2.8676E-01 & & 6.7967E-10 \\
		9.0000E-01 & & 1.5054E-01 & & 1.5054E-01 & & 7.2359E-10 \\
		9.5000E-01 & & 7.7154E-02 & & 7.7154E-02 & & 6.3819E-10 \\
		9.9000E-01 & & 1.5741E-02 & & 1.5741E-02 & & 1.2248E-10 \\
		1.0000E+00 & & 0.0000E+00 & & 0.0000E+00 & & 5.4216E-18 \\
		\hline\hline
	\end{tabular}
\end{table}

\begin{table}[H]
	\centering
	\scriptsize
	\caption{Pointwise LWM results for $\mathrm{Pe}=10$.}
	\label{tab:lwm_pe10_pointwise}
	\begin{tabular}{||lllllll||}
		\hline\hline
		$x$ & & LWM solution & & Exact solution & & Absolute error \\
		\hline\hline
		0.0000E+00 & & 1.0000E+00 & & 1.0000E+00 & & 0.0000E+00 \\
		1.0000E-01 & & 9.9992E-01 & & 9.9992E-01 & & 3.9500E-10 \\
		2.0000E-01 & & 9.9971E-01 & & 9.9971E-01 & & 3.1653E-09 \\
		3.0000E-01 & & 9.9913E-01 & & 9.9913E-01 & & 7.9017E-09 \\
		4.0000E-01 & & 9.9757E-01 & & 9.9757E-01 & & 2.8303E-08 \\
		5.0000E-01 & & 9.9331E-01 & & 9.9331E-01 & & 8.1627E-09 \\
		6.0000E-01 & & 9.8173E-01 & & 9.8173E-01 & & 6.4677E-08 \\
		7.0000E-01 & & 9.5026E-01 & & 9.5026E-01 & & 4.2507E-07 \\
		8.0000E-01 & & 8.6470E-01 & & 8.6470E-01 & & 1.0230E-06 \\
		9.0000E-01 & & 6.3215E-01 & & 6.3215E-01 & & 4.6358E-06 \\
		9.5000E-01 & & 3.9348E-01 & & 3.9349E-01 & & 4.8073E-06 \\
		9.8000E-01 & & 1.8128E-01 & & 1.8128E-01 & & 3.8032E-06 \\
		9.9000E-01 & & 9.5169E-02 & & 9.5167E-02 & & 2.1373E-06 \\
		9.9500E-01 & & 4.8774E-02 & & 4.8773E-02 & & 7.4286E-07 \\
		9.9900E-01 & & 9.9507E-03 & & 9.9506E-03 & & 3.7768E-08 \\
		1.0000E+00 & & 0.0000E+00 & & 0.0000E+00 & & 3.6253E-18 \\
		\hline\hline
	\end{tabular}
\end{table}

\begin{table}[H]
	\centering
	\scriptsize
	\caption{Pointwise LWM results for $\mathrm{Pe}=100$.}
	\label{tab:lwm_pe100_pointwise}
	\begin{tabular}{||lllllll||}
		\hline\hline
		$x$ & & LWM solution & & Exact solution & & Absolute error \\
		\hline\hline
		0.0000E+00 & & 1.0000E+00 & & 1.0000E+00 & & 0.0000E+00 \\
		1.0000E-01 & & 1.0000E+00 & & 1.0000E+00 & & 0.0000E+00 \\
		2.0000E-01 & & 1.0000E+00 & & 1.0000E+00 & & 0.0000E+00 \\
		3.0000E-01 & & 1.0000E+00 & & 1.0000E+00 & & 0.0000E+00 \\
		4.0000E-01 & & 1.0000E+00 & & 1.0000E+00 & & 0.0000E+00 \\
		5.0000E-01 & & 1.0000E+00 & & 1.0000E+00 & & 0.0000E+00 \\
		6.0000E-01 & & 1.0000E+00 & & 1.0000E+00 & & 0.0000E+00 \\
		7.0000E-01 & & 1.0000E+00 & & 1.0000E+00 & & 1.1102E-16 \\
		8.0000E-01 & & 1.0000E+00 & & 1.0000E+00 & & 2.9665E-12 \\
		9.0000E-01 & & 9.9995E-01 & & 9.9995E-01 & & 3.2658E-08 \\
		9.5000E-01 & & 9.9326E-01 & & 9.9326E-01 & & 2.4230E-06 \\
		9.8000E-01 & & 8.6485E-01 & & 8.6466E-01 & & 1.8280E-04 \\
		9.9000E-01 & & 6.3194E-01 & & 6.3212E-01 & & 1.8429E-04 \\
		9.9500E-01 & & 3.9358E-01 & & 3.9347E-01 & & 1.1011E-04 \\
		9.9800E-01 & & 1.8135E-01 & & 1.8127E-01 & & 7.6448E-05 \\
		9.9900E-01 & & 9.5190E-02 & & 9.5163E-02 & & 2.7124E-05 \\
		9.9950E-01 & & 4.8779E-02 & & 4.8771E-02 & & 7.9670E-06 \\
		9.9990E-01 & & 9.9505E-03 & & 9.9502E-03 & & 3.6053E-07 \\
		1.0000E+00 & & 0.0000E+00 & & 0.0000E+00 & & 1.7738E-16 \\
		\hline\hline
	\end{tabular}
\end{table}

\begin{table}[H]
	\centering
	\scriptsize
	\caption{Pointwise LWM results for $\mathrm{Pe}=1000$.}
	\label{tab:lwm_pe1000_pointwise}
	\begin{tabular}{||lllllll||}
		\hline\hline
		$x$ & & LWM solution & & Exact solution & & Absolute error \\
		\hline\hline
		0.0000E+00 & & 1.0000E+00 & & 1.0000E+00 & & 0.0000E+00 \\
		1.0000E-01 & & 1.0000E+00 & & 1.0000E+00 & & 1.1102E-15 \\
		2.0000E-01 & & 1.0000E+00 & & 1.0000E+00 & & 3.4417E-15 \\
		3.0000E-01 & & 1.0000E+00 & & 1.0000E+00 & & 5.3291E-15 \\
		4.0000E-01 & & 1.0000E+00 & & 1.0000E+00 & & 7.5495E-15 \\
		5.0000E-01 & & 1.0000E+00 & & 1.0000E+00 & & 9.5479E-15 \\
		6.0000E-01 & & 1.0000E+00 & & 1.0000E+00 & & 1.1879E-14 \\
		7.0000E-01 & & 1.0000E+00 & & 1.0000E+00 & & 1.3989E-14 \\
		8.0000E-01 & & 1.0000E+00 & & 1.0000E+00 & & 1.6209E-14 \\
		9.0000E-01 & & 1.0000E+00 & & 1.0000E+00 & & 1.8319E-14 \\
		9.9500E-01 & & 9.9223E-01 & & 9.9326E-01 & & 1.0370E-03 \\
		9.9800E-01 & & 8.6027E-01 & & 8.6466E-01 & & 4.3907E-03 \\
		9.9900E-01 & & 6.3414E-01 & & 6.3212E-01 & & 2.0149E-03 \\
		9.9950E-01 & & 3.9570E-01 & & 3.9347E-01 & & 2.2267E-03 \\
		9.9980E-01 & & 1.8192E-01 & & 1.8127E-01 & & 6.5245E-04 \\
		9.9990E-01 & & 9.5358E-02 & & 9.5163E-02 & & 1.9511E-04 \\
		9.9995E-01 & & 4.8824E-02 & & 4.8771E-02 & & 5.3173E-05 \\
		9.9999E-01 & & 9.9524E-03 & & 9.9502E-03 & & 2.2756E-06 \\
		1.0000E+00 & & 0.0000E+00 & & 0.0000E+00 & & 7.8063E-17 \\
		\hline\hline
	\end{tabular}
\end{table}

Table~\ref{tab:lwm_summary} reports the main numerical summary for the direct two-sum Legendre wavelet method. The algebraic residual $\|A\mathbf{c}-\mathbf{b}\|_{\infty}$ remains small in all four cases, showing that the assembled linear systems are solved accurately. The condition number increases with $\mathrm{Pe}$, indicating stronger ill-conditioning of the algebraic system in the boundary-layer regime.

\begin{table}[H]
	\centering
	\scriptsize
	\caption{Summary of direct two-sum Legendre wavelet results.}
	\label{tab:lwm_summary}
	\begin{tabular}{||lllllllllll||}
		\hline\hline
		Pe & $\nu$ & $J$ & $M$ & Unknowns & cond$(A)$ & Alg. Res. & $E_{\infty}$ & $E_{\mathrm{mean}}$ & $E_{\mathrm{rms}}$ & $E_2$ \\
		\hline\hline
		1    & 1.000E+00 & 4   & 6 & 24   & 2.949E+04 & 4.441E-16 & 1.056E-09 & 1.485E-10 & 2.821E-10 & 4.341E-10 \\
		10   & 1.000E-01 & 8   & 6 & 48   & 3.243E+05 & 2.842E-14 & 6.327E-06 & 3.779E-07 & 1.117E-06 & 1.312E-06 \\
		100  & 1.000E-02 & 40  & 6 & 240  & 4.458E+07 & 1.819E-12 & 2.213E-04 & 8.397E-06 & 3.218E-05 & 2.021E-05 \\
		1000 & 1.000E-03 & 200 & 6 & 1200 & 5.879E+09 & 1.310E-10 & 4.735E-03 & 1.708E-04 & 6.758E-04 & 2.006E-04 \\
		\hline\hline
	\end{tabular}
\end{table}

The LWM results show that the direct two-sum Legendre wavelet method captures the solution profile for all four Peclet numbers. For $\mathrm{Pe}=1$ and $\mathrm{Pe}=10$, the errors remain small throughout the interval. For $\mathrm{Pe}=100$ and $\mathrm{Pe}=1000$, the dominant errors are concentrated near the boundary layer close to $x=1$. Although the maximum error increases with $\mathrm{Pe}$, the method still resolves the sharp transition for the strongest case considered here, with a maximum absolute error of order $10^{-3}$.

\subsection{Standard PINN results}

The standard soft-boundary PINN is tested on the same four Peclet numbers. The network uses four hidden layers, 64 neurons in each hidden layer, and the hyperbolic tangent activation function. The residual loss is evaluated using $10000$ uniformly distributed collocation points during training, and the trained network is then evaluated on the same dense grid used for the LWM results.

\begin{figure}[H]
	\centering
	\begin{minipage}{0.48\textwidth}
		\centering
		\includegraphics[width=\textwidth]{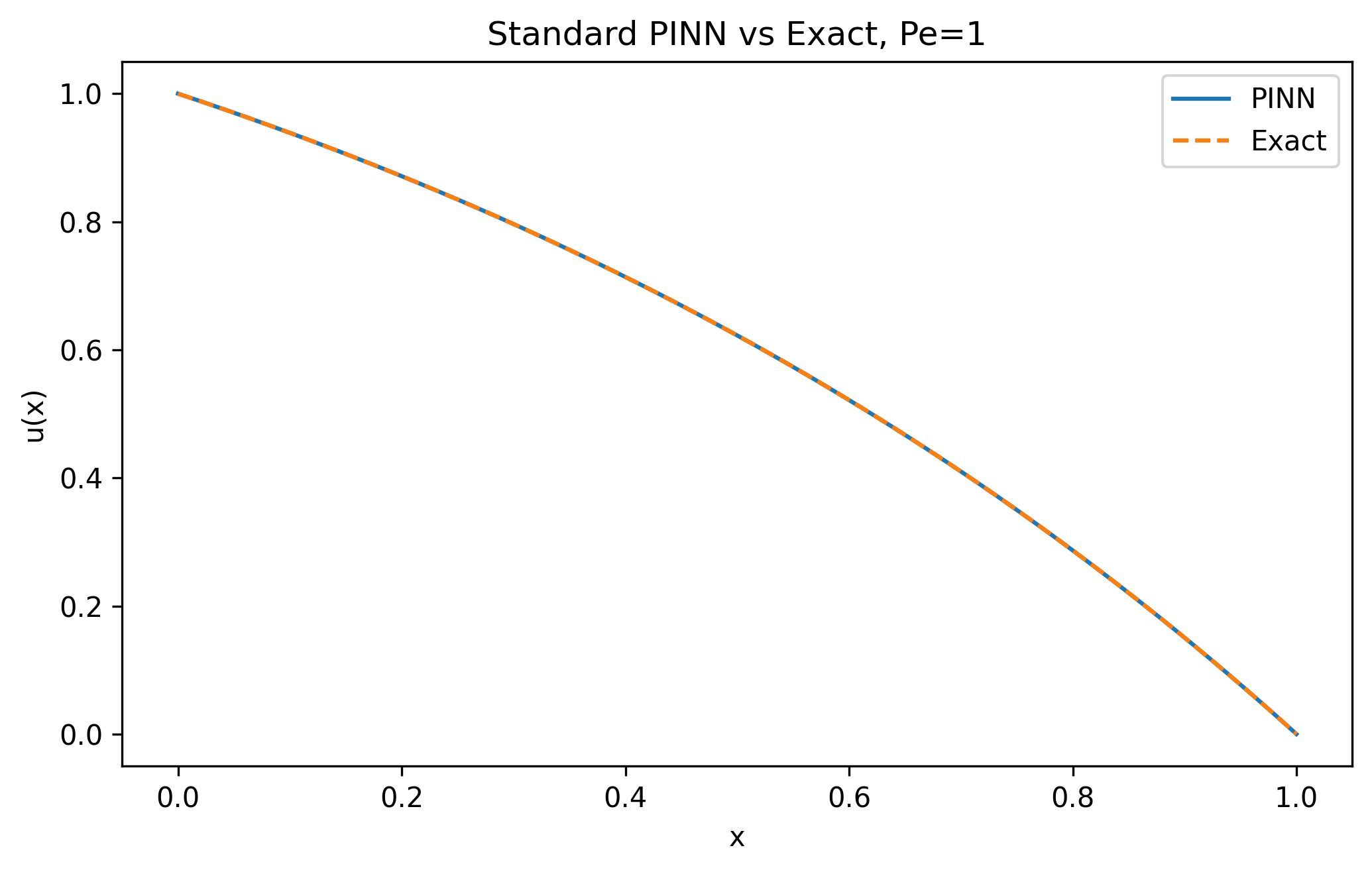}\\
		\textbf{(a)} $\mathrm{Pe}=1$
	\end{minipage}
	\begin{minipage}{0.48\textwidth}
		\centering
		\includegraphics[width=\textwidth]{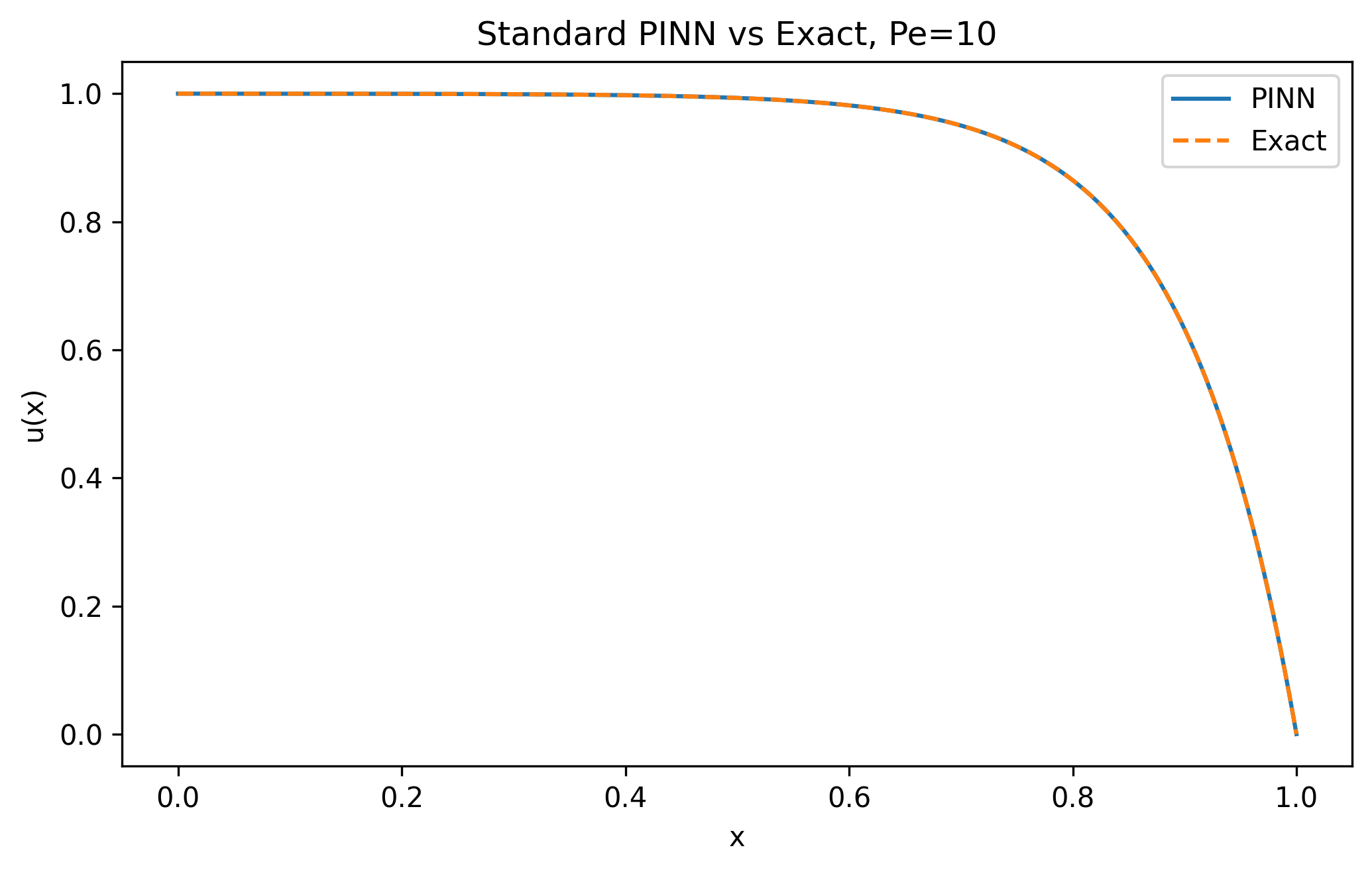}\\
		\textbf{(b)} $\mathrm{Pe}=10$
	\end{minipage}
	
	\vspace{2mm}
	
	\begin{minipage}{0.48\textwidth}
		\centering
		\includegraphics[width=\textwidth]{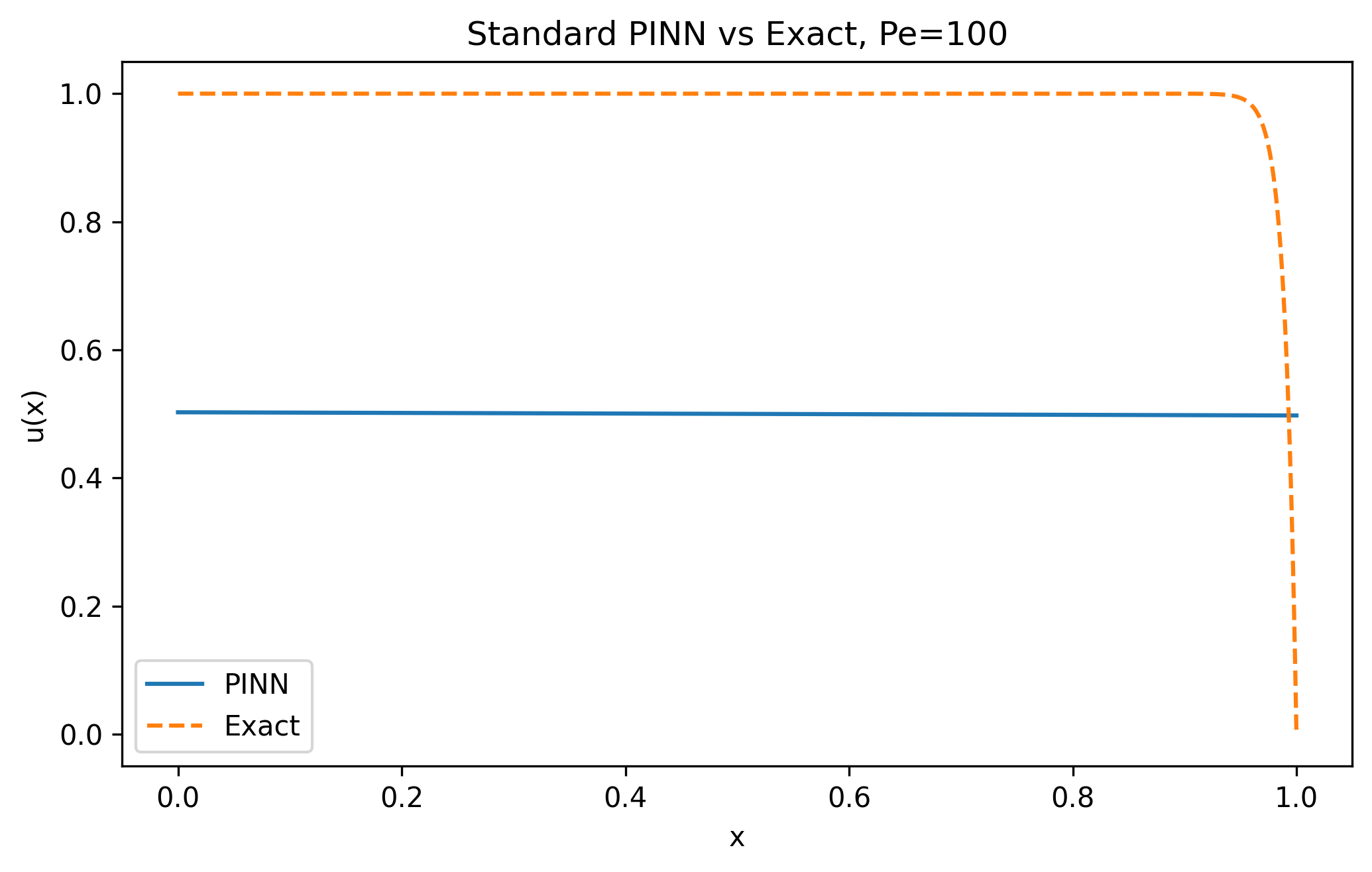}\\
		\textbf{(c)} $\mathrm{Pe}=100$
	\end{minipage}
	\begin{minipage}{0.48\textwidth}
		\centering
		\includegraphics[width=\textwidth]{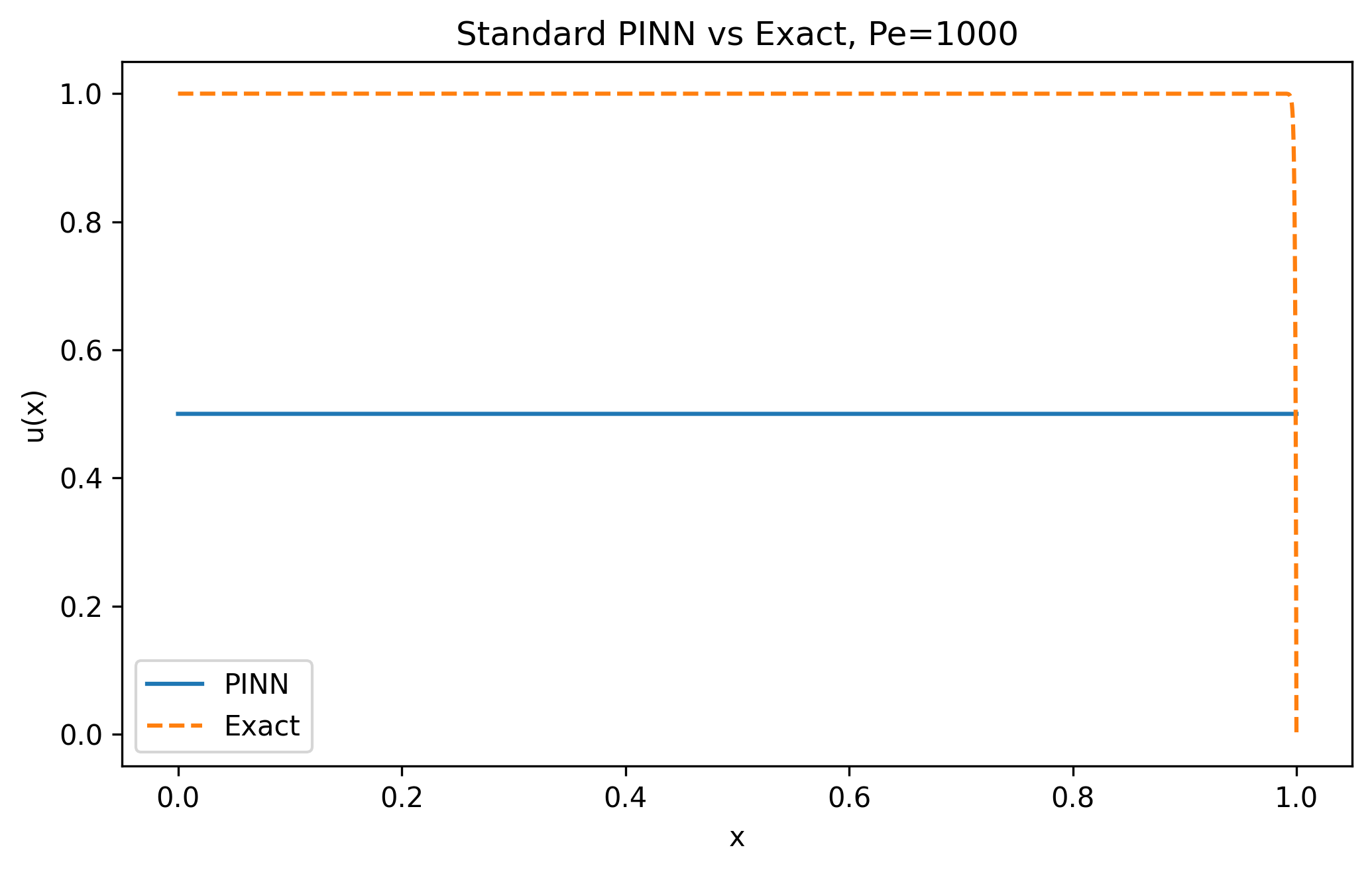}\\
		\textbf{(d)} $\mathrm{Pe}=1000$
	\end{minipage}
	\caption{Standard PINN solution and exact solution for different Peclet numbers.}
	\label{fig:pinn_solution_all}
\end{figure}

Figure~\ref{fig:pinn_solution_all} shows the standard PINN solution together with the exact solution. For $\mathrm{Pe}=1$ and $\mathrm{Pe}=10$, the PINN approximation follows the exact solution closely. For $\mathrm{Pe}=100$ and $\mathrm{Pe}=1000$, the approximation becomes nearly constant and does not resolve the sharp boundary layer near $x=1$.

\begin{figure}[H]
	\centering
	\begin{minipage}{0.48\textwidth}
		\centering
		\includegraphics[width=\textwidth]{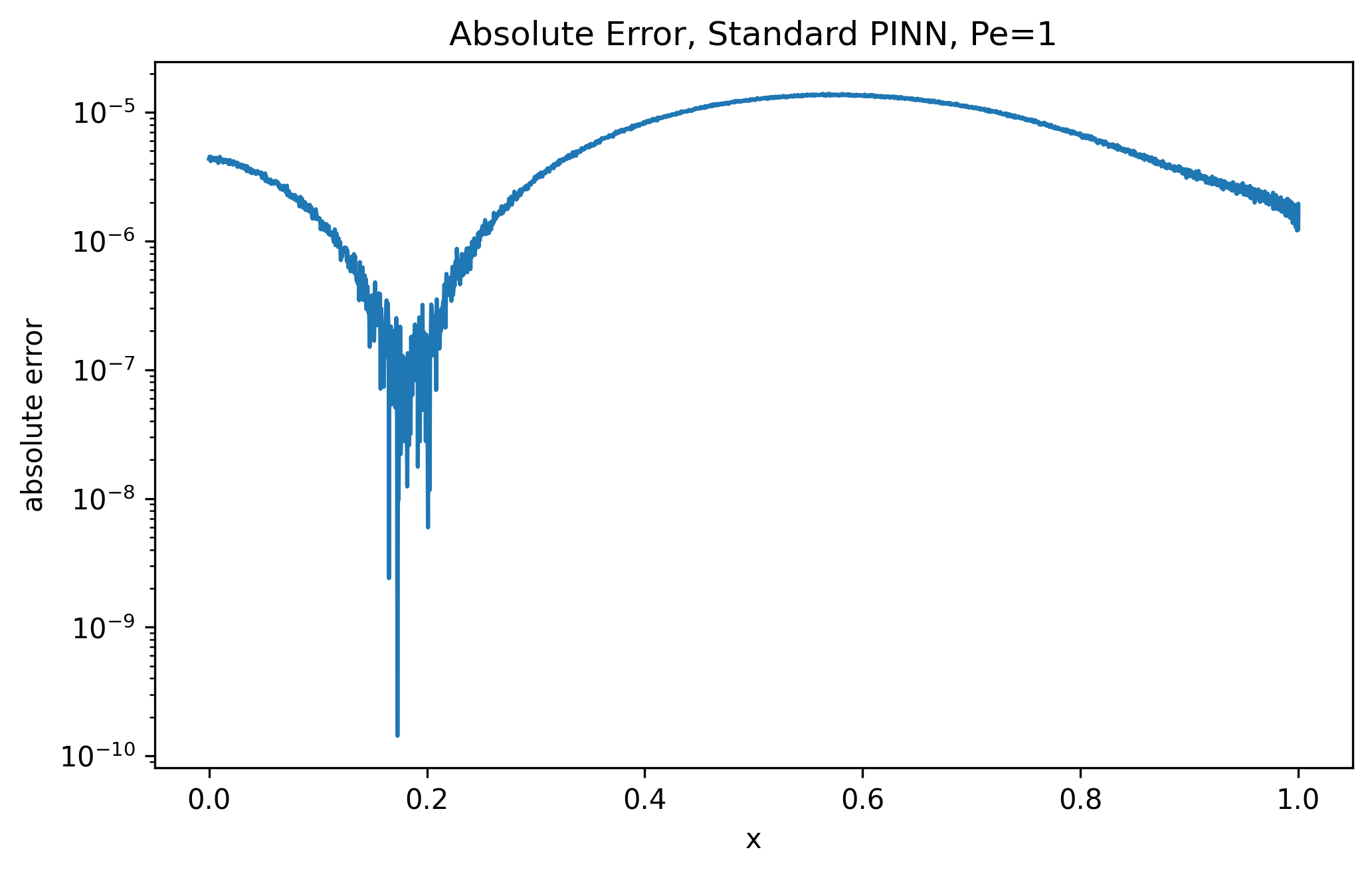}\\
		\textbf{(a)} $\mathrm{Pe}=1$
	\end{minipage}
	\begin{minipage}{0.48\textwidth}
		\centering
		\includegraphics[width=\textwidth]{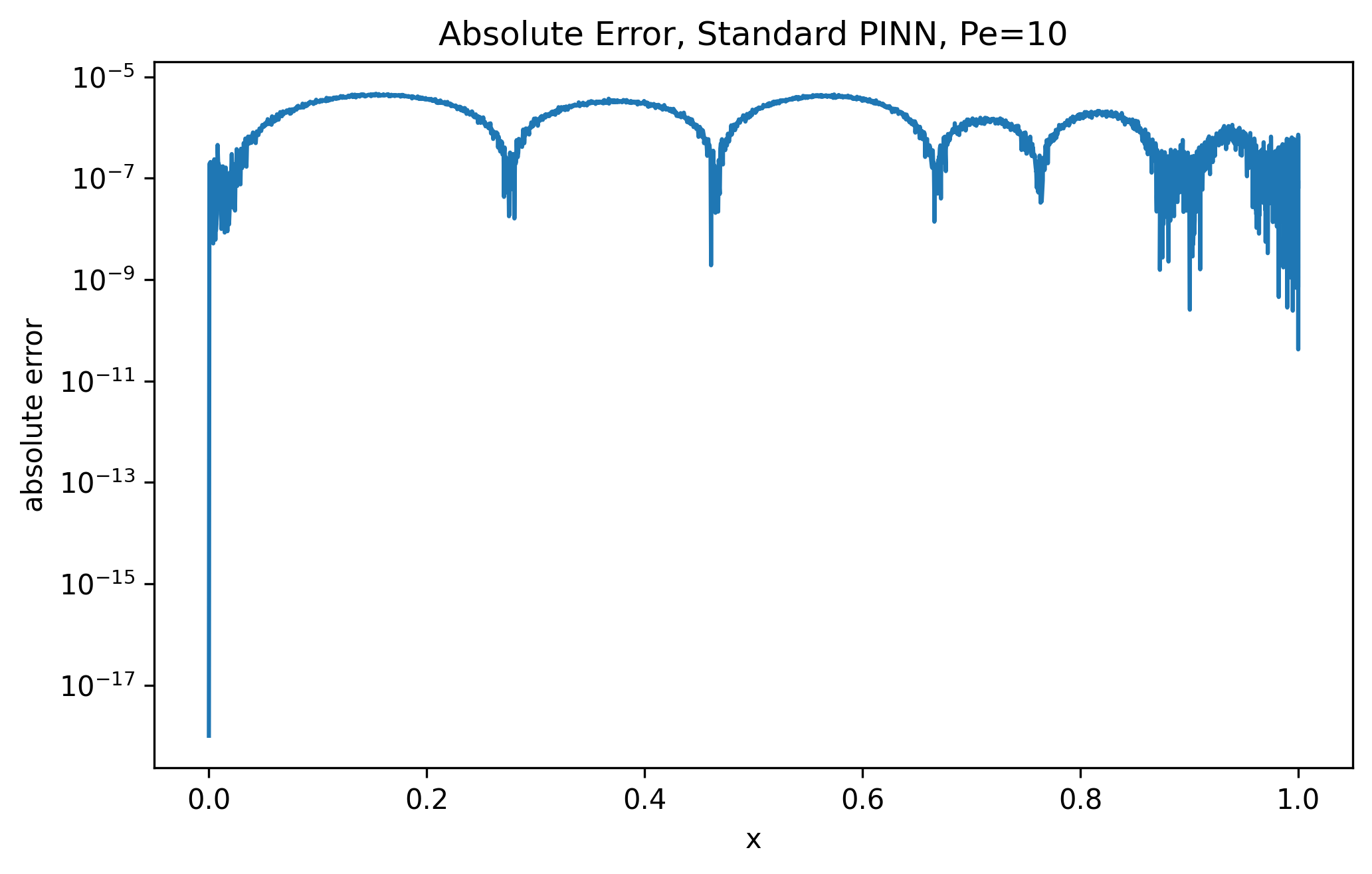}\\
		\textbf{(b)} $\mathrm{Pe}=10$
	\end{minipage}
	
	\vspace{2mm}
	
	\begin{minipage}{0.48\textwidth}
		\centering
		\includegraphics[width=\textwidth]{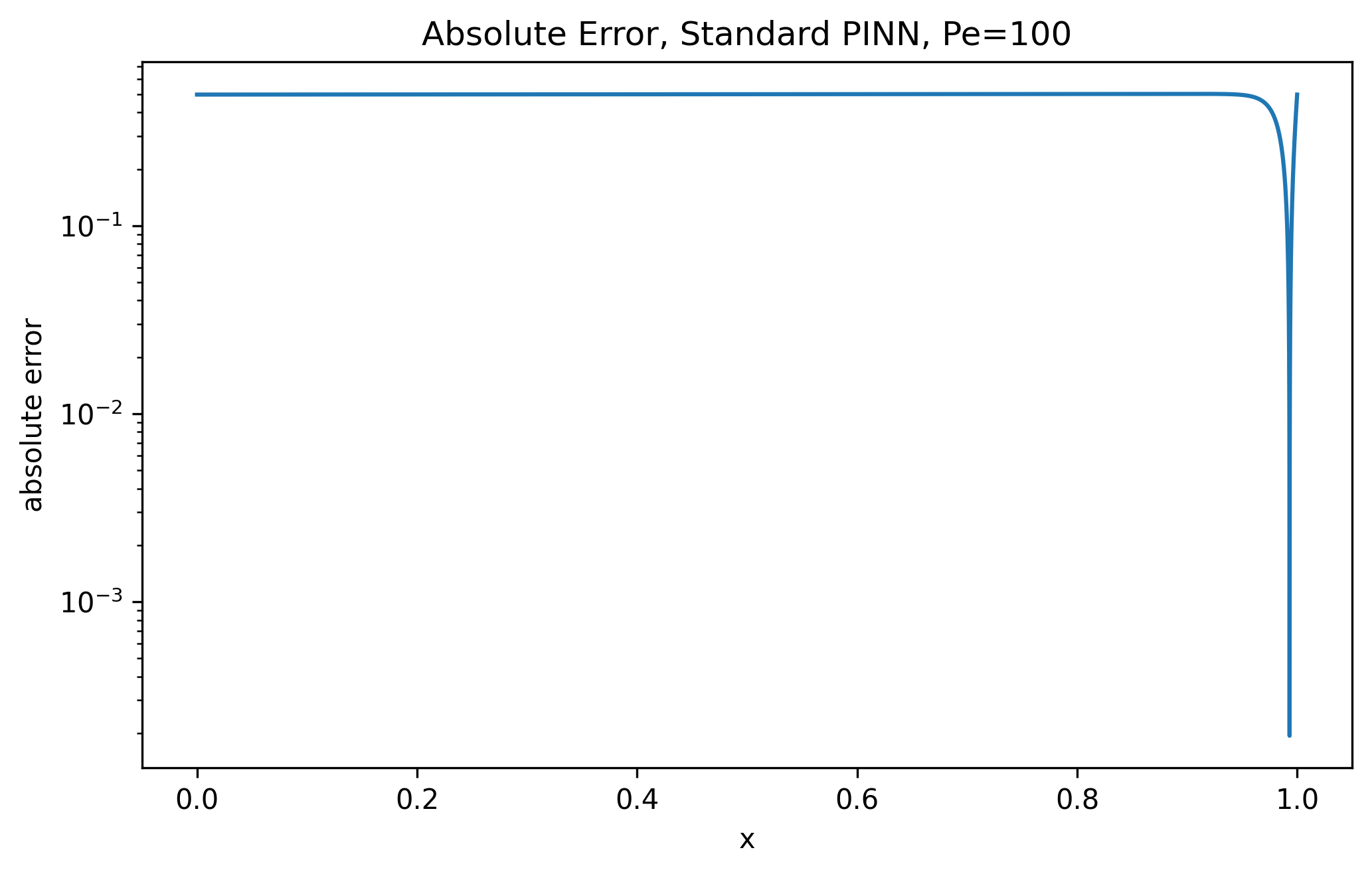}\\
		\textbf{(c)} $\mathrm{Pe}=100$
	\end{minipage}
	\begin{minipage}{0.48\textwidth}
		\centering
		\includegraphics[width=\textwidth]{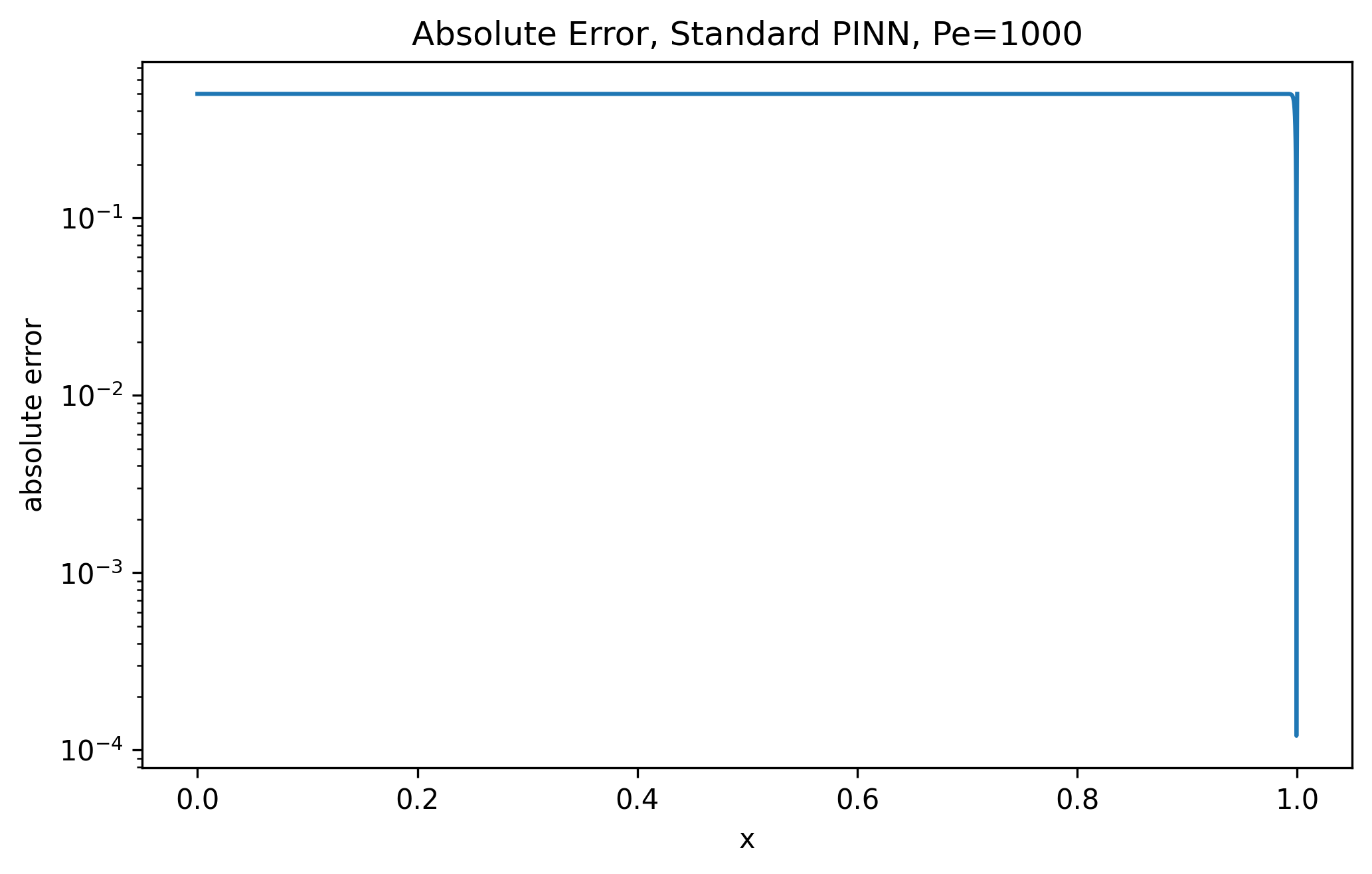}\\
		\textbf{(d)} $\mathrm{Pe}=1000$
	\end{minipage}
	\caption{Absolute error of the standard PINN for different Peclet numbers.}
	\label{fig:pinn_error_all}
\end{figure}

Figure~\ref{fig:pinn_error_all} shows the absolute-error profiles for the standard PINN. The error remains small for $\mathrm{Pe}=1$ and $\mathrm{Pe}=10$. For $\mathrm{Pe}=100$ and $\mathrm{Pe}=1000$, the error is much larger because the trained network does not capture the boundary-layer profile.

The pointwise PINN results are reported in Tables~\ref{tab:pinn_pe1_pointwise}--\ref{tab:pinn_pe1000_pointwise}. The tabulated points follow the same near-boundary evaluation pattern used in the LWM tables.

\begin{table}[H]
	\centering
	\scriptsize
	\caption{Pointwise PINN results for $\mathrm{Pe}=1$.}
	\label{tab:pinn_pe1_pointwise}
	\begin{tabular}{||lllllll||}
		\hline\hline
		$x$ & & PINN solution & & Exact solution & & Absolute error \\
		\hline\hline
		0.0000E+00 & & 1.0000E+00 & & 1.0000E+00 & & 4.5896E-06 \\
		1.0000E-01 & & 9.3879E-01 & & 9.3879E-01 & & 1.5216E-06 \\
		2.0000E-01 & & 8.7115E-01 & & 8.7115E-01 & & 1.0604E-07 \\
		3.0000E-01 & & 7.9639E-01 & & 7.9639E-01 & & 3.0681E-06 \\
		4.0000E-01 & & 7.1376E-01 & & 7.1377E-01 & & 8.2121E-06 \\
		5.0000E-01 & & 6.2245E-01 & & 6.2246E-01 & & 1.2615E-05 \\
		6.0000E-01 & & 5.2153E-01 & & 5.2155E-01 & & 1.3532E-05 \\
		7.0000E-01 & & 4.1001E-01 & & 4.1002E-01 & & 1.0928E-05 \\
		8.0000E-01 & & 2.8676E-01 & & 2.8676E-01 & & 6.6744E-06 \\
		9.0000E-01 & & 1.5054E-01 & & 1.5054E-01 & & 3.3994E-06 \\
		9.5000E-01 & & 7.7151E-02 & & 7.7154E-02 & & 2.5116E-06 \\
		9.9000E-01 & & 1.5739E-02 & & 1.5741E-02 & & 1.7436E-06 \\
		1.0000E+00 & & -1.5795E-06 & & 0.0000E+00 & & 1.5795E-06 \\
		\hline\hline
	\end{tabular}
\end{table}

\begin{table}[H]
	\centering
	\scriptsize
	\caption{Pointwise PINN results for $\mathrm{Pe}=10$.}
	\label{tab:pinn_pe10_pointwise}
	\begin{tabular}{||lllllll||}
		\hline\hline
		$x$ & & PINN solution & & Exact solution & & Absolute error \\
		\hline\hline
		0.0000E+00 & & 1.0000E+00 & & 1.0000E+00 & & 2.9802E-07 \\
		1.0000E-01 & & 9.9993E-01 & & 9.9992E-01 & & 3.2096E-06 \\
		2.0000E-01 & & 9.9971E-01 & & 9.9971E-01 & & 3.7948E-06 \\
		3.0000E-01 & & 9.9913E-01 & & 9.9913E-01 & & 1.2626E-06 \\
		4.0000E-01 & & 9.9756E-01 & & 9.9757E-01 & & 3.1155E-06 \\
		5.0000E-01 & & 9.9331E-01 & & 9.9331E-01 & & 2.0507E-06 \\
		6.0000E-01 & & 9.8173E-01 & & 9.8173E-01 & & 3.4965E-06 \\
		7.0000E-01 & & 9.5025E-01 & & 9.5026E-01 & & 1.3349E-06 \\
		8.0000E-01 & & 8.6471E-01 & & 8.6470E-01 & & 1.7671E-06 \\
		9.0000E-01 & & 6.3215E-01 & & 6.3215E-01 & & 1.5806E-07 \\
		9.5000E-01 & & 3.9349E-01 & & 3.9349E-01 & & 6.6611E-07 \\
		9.8000E-01 & & 1.8128E-01 & & 1.8128E-01 & & 3.2103E-07 \\
		9.9000E-01 & & 9.5167E-02 & & 9.5167E-02 & & 4.7974E-08 \\
		9.9500E-01 & & 4.8773E-02 & & 4.8773E-02 & & 1.4659E-08 \\
		9.9900E-01 & & 9.9505E-03 & & 9.9506E-03 & & 1.6646E-07 \\
		1.0000E+00 & & -1.2666E-07 & & 0.0000E+00 & & 1.2666E-07 \\
		\hline\hline
	\end{tabular}
\end{table}

\begin{table}[H]
	\centering
	\scriptsize
	\caption{Pointwise PINN results for $\mathrm{Pe}=100$.}
	\label{tab:pinn_pe100_pointwise}
	\begin{tabular}{||lllllll||}
		\hline\hline
		$x$ & & PINN solution & & Exact solution & & Absolute error \\
		\hline\hline
		0.0000E+00 & & 5.0247E-01 & & 1.0000E+00 & & 4.9753E-01 \\
		1.0000E-01 & & 5.0200E-01 & & 1.0000E+00 & & 4.9800E-01 \\
		2.0000E-01 & & 5.0152E-01 & & 1.0000E+00 & & 4.9848E-01 \\
		3.0000E-01 & & 5.0104E-01 & & 1.0000E+00 & & 4.9896E-01 \\
		4.0000E-01 & & 5.0055E-01 & & 1.0000E+00 & & 4.9945E-01 \\
		5.0000E-01 & & 5.0005E-01 & & 1.0000E+00 & & 4.9995E-01 \\
		6.0000E-01 & & 4.9955E-01 & & 1.0000E+00 & & 5.0045E-01 \\
		7.0000E-01 & & 4.9905E-01 & & 1.0000E+00 & & 5.0095E-01 \\
		8.0000E-01 & & 4.9854E-01 & & 1.0000E+00 & & 5.0146E-01 \\
		9.0000E-01 & & 4.9802E-01 & & 9.9995E-01 & & 5.0194E-01 \\
		9.5000E-01 & & 4.9776E-01 & & 9.9326E-01 & & 4.9551E-01 \\
		9.8000E-01 & & 4.9760E-01 & & 8.6466E-01 & & 3.6707E-01 \\
		9.9000E-01 & & 4.9754E-01 & & 6.3212E-01 & & 1.3458E-01 \\
		9.9500E-01 & & 4.9752E-01 & & 3.9347E-01 & & 1.0405E-01 \\
		9.9800E-01 & & 4.9750E-01 & & 1.8127E-01 & & 3.1623E-01 \\
		9.9900E-01 & & 4.9750E-01 & & 9.5163E-02 & & 4.0233E-01 \\
		9.9950E-01 & & 4.9749E-01 & & 4.8771E-02 & & 4.4872E-01 \\
		9.9990E-01 & & 4.9749E-01 & & 9.9502E-03 & & 4.8754E-01 \\
		1.0000E+00 & & 4.9749E-01 & & 0.0000E+00 & & 4.9749E-01 \\
		\hline\hline
	\end{tabular}
\end{table}

\begin{table}[H]
	\centering
	\scriptsize
	\caption{Pointwise PINN results for $\mathrm{Pe}=1000$.}
	\label{tab:pinn_pe1000_pointwise}
	\begin{tabular}{||lllllll||}
		\hline\hline
		$x$ & & PINN solution & & Exact solution & & Absolute error \\
		\hline\hline
		0.0000E+00 & & 5.0004E-01 & & 1.0000E+00 & & 4.9996E-01 \\
		1.0000E-01 & & 5.0004E-01 & & 1.0000E+00 & & 4.9996E-01 \\
		2.0000E-01 & & 5.0003E-01 & & 1.0000E+00 & & 4.9997E-01 \\
		3.0000E-01 & & 5.0003E-01 & & 1.0000E+00 & & 4.9997E-01 \\
		4.0000E-01 & & 5.0002E-01 & & 1.0000E+00 & & 4.9998E-01 \\
		5.0000E-01 & & 5.0002E-01 & & 1.0000E+00 & & 4.9998E-01 \\
		6.0000E-01 & & 5.0001E-01 & & 1.0000E+00 & & 4.9999E-01 \\
		7.0000E-01 & & 5.0001E-01 & & 1.0000E+00 & & 4.9999E-01 \\
		8.0000E-01 & & 5.0000E-01 & & 1.0000E+00 & & 5.0000E-01 \\
		9.0000E-01 & & 5.0000E-01 & & 1.0000E+00 & & 5.0000E-01 \\
		9.9500E-01 & & 4.9999E-01 & & 9.9326E-01 & & 4.9327E-01 \\
		9.9800E-01 & & 4.9999E-01 & & 8.6466E-01 & & 3.6467E-01 \\
		9.9900E-01 & & 4.9999E-01 & & 6.3212E-01 & & 1.3213E-01 \\
		9.9950E-01 & & 4.9999E-01 & & 3.9347E-01 & & 1.0652E-01 \\
		9.9980E-01 & & 4.9999E-01 & & 1.8127E-01 & & 3.1872E-01 \\
		9.9990E-01 & & 4.9999E-01 & & 9.5163E-02 & & 4.0483E-01 \\
		9.9995E-01 & & 4.9999E-01 & & 4.8771E-02 & & 4.5122E-01 \\
		9.9999E-01 & & 4.9999E-01 & & 9.9502E-03 & & 4.9004E-01 \\
		1.0000E+00 & & 4.9999E-01 & & 0.0000E+00 & & 4.9999E-01 \\
		\hline\hline
	\end{tabular}
\end{table}

Table~\ref{tab:pinn_summary} reports the summary of the standard PINN results. The standard formulation gives small errors for $\mathrm{Pe}=1$ and $\mathrm{Pe}=10$. However, for $\mathrm{Pe}=100$ and $\mathrm{Pe}=1000$, the maximum absolute error is approximately $5.0\mathrm{E}{-01}$, which indicates that the sharp boundary layer is not resolved.

\begin{table}[H]
	\centering
	\scriptsize
	\caption{Summary of standard PINN results.}
	\label{tab:pinn_summary}
	\begin{tabular}{||lllllllllll||}
		\hline\hline
		Pe & $\nu$ & Total loss & Res. loss & BC loss & $E_{\infty}$ & $E_{\mathrm{mean}}$ & $E_{\mathrm{rms}}$ & $E_2$ & Left BE & Right BE \\
		\hline\hline
		1    & 1.000E+00 & 8.763E-08 & 8.541E-08 & 2.225E-11 & 1.391E-05 & 3.235E-06 & 4.760E-06 & 7.792E-06 & 4.470E-06 & 1.505E-06 \\
		10   & 1.000E-01 & 1.145E-06 & 1.145E-06 & 3.469E-14 & 4.592E-06 & 7.156E-07 & 1.332E-06 & 2.313E-06 & 1.788E-07 & 5.215E-08 \\
		100  & 1.000E-02 & 4.975E+01 & 2.477E-01 & 4.950E-01 & 5.019E-01 & 4.751E-01 & 4.815E-01 & 4.949E-01 & 4.975E-01 & 4.975E-01 \\
		1000 & 1.000E-03 & 5.000E+01 & 2.603E-03 & 4.999E-01 & 5.000E-01 & 4.783E-01 & 4.842E-01 & 4.995E-01 & 5.000E-01 & 5.000E-01 \\
		\hline\hline
	\end{tabular}
\end{table}

\subsection{Direct comparison between LWM and standard PINN}

Table~\ref{tab:lwm_pinn_comparison} compares the main error measures of the direct two-sum Legendre wavelet method and the standard PINN. For $\mathrm{Pe}=1$, both methods approximate the solution, with the LWM giving the smaller maximum error. For $\mathrm{Pe}=10$, both methods give small errors. For $\mathrm{Pe}=100$ and $\mathrm{Pe}=1000$, the difference becomes substantial: the LWM retains much smaller errors, while the standard PINN gives maximum errors of approximately $5.0\times 10^{-1}$.

\begin{table}[H]
	\centering
	\scriptsize
	\caption{Direct comparison of LWM and standard PINN errors.}
	\label{tab:lwm_pinn_comparison}
    \begin{tabular}{||llllll||}
		\hline\hline
		Pe & $\nu$ & LWM $E_{\infty}$ & PINN $E_{\infty}$ & LWM $E_2$ & PINN $E_2$ \\
		\hline\hline
		1    & 1.000E+00 & 1.056E-09 & 1.391E-05 & 4.341E-10 & 7.792E-06 \\
		10   & 1.000E-01 & 6.327E-06 & 4.592E-06 & 1.312E-06 & 2.313E-06 \\
		100  & 1.000E-02 & 2.213E-04 & 5.019E-01 & 2.021E-05 & 4.949E-01 \\
		1000 & 1.000E-03 & 4.735E-03 & 5.000E-01 & 2.006E-04 & 4.995E-01 \\
		\hline\hline
	\end{tabular}
\end{table}

\subsection{Discussion}

The results indicate that the difference between the two methods is mainly caused by how each approximation represents localized structure. The Legendre wavelet method uses locally supported polynomial basis functions, so increasing the number of cells directly increases resolution near the boundary layer. This gives the method a simple refinement mechanism for the increasingly thin transition region as $\mathrm{Pe}$ increases.

The standard soft-boundary PINN does not have this local refinement mechanism. In the high-Peclet cases, the loss minimization leads to a nearly constant approximation instead of the sharp transition near $x=1$. This behavior is a limitation of the baseline formulation tested here, where the boundary conditions are imposed only through penalty terms and no hard constraint, adaptive weighting, domain decomposition, or boundary-layer correction is used.

The comparison also shows that evaluation strategy matters for singularly perturbed problems. Since the main error is concentrated in a narrow region near the right boundary, dense evaluation with additional near-boundary points is necessary. Without such evaluation, the apparent accuracy can be misleading because most of the interval is almost constant for large $\mathrm{Pe}$.

\section{Conclusion}

This work compared a direct two-sum Legendre wavelet collocation method and a standard soft-boundary PINN for a singularly perturbed advection--diffusion boundary-value problem. The results show that the Legendre wavelet method remains effective across the tested Peclet numbers, while the standard soft-boundary PINN loses accuracy in the sharper boundary-layer cases.

Overall, the study demonstrates that local algebraic collocation provides a more dependable approximation mechanism than the baseline residual-loss PINN formulation for this boundary-layer benchmark. Future work will test hard-boundary, adaptive-collocation, and domain-decomposed PINN variants on the same problem.

\subsection*{Code availability}

The numerical experiments were implemented in Python using NumPy for the Legendre wavelet collocation method and PyTorch for the PINN computations. The source codes and numerical data will be made available in a public GitHub repository.

	\end{document}